\documentclass[preprint,12pt,3p]{elsarticle}

\usepackage{inputenc}
\usepackage{amsmath,amssymb,amsfonts,color,bm}
\usepackage{graphicx}
\usepackage{epstopdf}
\usepackage{latexsym}
\usepackage{multirow}
\usepackage{float}
\usepackage{verbatim}
\usepackage{enumitem}
\usepackage{algorithmic}
\usepackage{amsopn}
\usepackage{enumitem}
\usepackage{tikz-cd}
\usepackage{comment}
\usepackage[scr=boondoxo]{mathalpha}
\usepackage{hyperref}
\usepackage{orcidlink}

\def\emm{\widehat{\ell}}
\def\R{{\mathbb R}}

\def\N{{\mathbb N}}

\def\P{{\mathbb P}}
\def\safek{\mathcal{S}^{\mathrm{Fek}}}

\newcommand{\de}{\mathrm{d}}
\newcommand{\op}{\mathrm{op}}

\ifpdf
  \DeclareGraphicsExtensions{.eps,.pdf,.png,.jpg}
\else
  \DeclareGraphicsExtensions{.eps}
\fi

\newtheorem{theorem}{Theorem}
\newtheorem{proposition}[theorem]{Proposition}
\newtheorem{remark}[theorem]{Remark}
\newtheorem{lemma}[theorem]{Lemma}
\newtheorem{problem}[theorem]{Problem}
\newtheorem{definition}[theorem]{Definition}

\newcommand*{\QED}{\null\nobreak\hfill\ensuremath{\square}\vspace{1mm}}%

\usepackage{xpatch}

\xpatchcmd{\proof}{\itshape}{\prooflabelfont}{}{}
\newcommand{\prooflabelfont}{\bfseries}

\begin{document}

\begin{frontmatter}



\title{The Fekete problem in segmental polynomial interpolation}

\author[inst1,inst2]{Ludovico Bruni Bruno \orcidlink{0000-0002-5246-8049}}

\affiliation[inst1]{organization={Dipartimento di Matematica \textquotedblleft Tullio Levi-Civita\textquotedblright, Università di Padova},%Department and Organization
            addressline={via Trieste, 63}, 
            city={Padova},
            postcode={35131},
            country={Italia}}

\author[inst1]{Wolfgang Erb \orcidlink{0000-0003-3541-5401}}

\affiliation[inst2]{organization={Istituto Nazionale di Alta Matematica \textquotedblleft Francesco Severi\textquotedblright},%Department and Organization
            addressline={Piazzale Aldo Moro, 5}, 
            city={Roma},
            postcode={00185}, 
            country={Italia}}

\begin{abstract}
In this article, we study the Fekete problem in segmental and combined nodal-segmental univariate polynomial interpolation by investigating sets of segments, or segments combined with nodes, such that the Vandermonde determinant for the respective polynomial interpolation problem is maximized. For particular families of segments, we will be able to find explicit solutions of the corresponding maximization problem. The quality of the Fekete segments depends hereby strongly on the utilized normalization of the segmental information in the Vandermonde matrix. To measure the quality of the Fekete segments in interpolation, we analyse the asymptotic behaviour of the generalized Lebesgue constant linked to the interpolation problem. For particular sets of Fekete segments we will get, similar to the nodal case, a favourable logarithmic growth of this constant.  
%We study and compare three interpolators and their features. Their norms can be characterised as three different Lebesgue constants: two of them are diffused in literature, whereas the remaining can be seen as a generalisation of the others. To each of these interpolators corresponds a different Fekete problem: two provide a concept of quasi-optimality of the set of supports, whereas the remaining one furnishes a negative answer to the following question: is it true that good interpolation nodes can be used to construct good interpolation simplices?
\end{abstract}

\begin{keyword}
Fekete problem \sep segmental polynomial interpolation \sep Lebesgue constant  \sep optimal design of segments \sep histopolation \sep polynomial approximation of differential forms
\MSC 41A05 \sep 41A25 \sep A1A30 \sep 65D05
\end{keyword}

\end{frontmatter}

\section{Introduction}

The quality of a polynomial interpolation operator depends strongly on the positioning of the nodes or the regions (referred to as supports) at which function samples or function averages are taken. The Runge phenomenon \cite{Runge} is a striking example of this: if the interpolation nodes are not wisely selected, the corresponding interpolating polynomial diverges sensibly from the interpolated function. For function averages over edge segments, similar effects have been encountered in \cite{BruniRunge}. Related phenomena have been observed also in very general situations \cite{BruniThesis}, evidencing the relevance of the geometrical structure and the distribution of the supports. 

Obviously, different supports determine different polynomial interpolators. To obtain an analytical measure for the quality of the geometry and the positioning of the chosen supports, we can use the operator norm of the interpolation operator which provides a general quantity for the numerical conditioning of the interpolation problem. This quantity is known as Lebesgue constant (see \cite{Ibrahimoglu} and the references therein) and is intimately related to the cardinal (or Lagrange) basis functions of the interpolation process. Since these basis functions may be obtained by Cramer's rule from the Vandermonde matrix, one deduces that supports that maximize the determinant of such a matrix provide a valuable set of supports also for interpolation. This is known as a \emph{Fekete problem} and has attracted a considerable amount of research in the case of nodal polynomial interpolation in the last century: the first works \cite{Fejer1932,Fejer1932b,Fekete} and recent works on Fekete problems \cite{BosSurvey,Bos2010,Bos1D,Bos2001,Hesthaven,Sommariva2009}, for example, span a period of about one hundred years. 

The generality of this reasoning raises the natural question whether this idea is applicable beyond the framework of pure nodal polynomial interpolation, for instance, to applications where one deals with averaged data \cite{Hiemstra2014}. This is particularly interesting in the one-dimensional setting, where one may replace nodal evaluations with integrals of a certain physical quantity along line segments \cite{BruniErb,Gerritsma} or by combined nodal and segmental data \cite{Robidoux2006}.  

In this work, we set up and study Fekete problems for these more general situations where segmental function averages or combined nodal-segmental data is used as input information. This results in the analysis and the comparison of three different polynomial interpolation operators: nodal, segmental, and combined nodal-segmental interpolation. The Fekete nodes in nodal polynomial interpolation are well-studied and will be used as a baseline for the comparison with the latter two scenarios. Segmental interpolation is also known as histopolation in the literature \cite{Gerritsma,Robidoux2006} and allows to consider interpolation also for less regular functions, whereas the third combined approach mixes the peculiarities of the first two. We dedicate Section \ref{sect:interpolators} to a brief discussion of the properties of these three families of interpolators. In Section \ref{sect:lebconst}, we review the features of the Lebesgue constants associated with the interpolators and show that there exists a common theory and terminology linking the three cases. Then, Section \ref{sect:feketeprobs} is devoted to the formalisation and resolution of the respective Fekete problems. While the nodal Fekete problem is classical and offers a wide literature, the segmental and the combined nodal-segmental Fekete problems have not been considered before. As not all degrees of freedom in the segmental Fekete maximization problem can be treated at once in an analytic way, we will restrict our attention to Fekete subproblems linked to particular families of supports as, for instance, families of concatenated segments or segments with uniform arc-length.  

The obtained results offer new and interesting insights in polynomial interpolation: while the nodal Fekete problem offers a quasi-optimal Lebesgue constant, the asymptotic behavior of the Lebesgue constant in the segmental case depends strongly on the chosen normalizations in the Vandermonde matrix. We will analyze this exemplarily in the case of concatenated segments. The obtained results further indicate that a commonly adopted technique for the construction of higher-dimensional supports, as for instance in \cite{BruniEdges}, which consists in the usage of good interpolation nodes as the endpoints of segments and simplices, leads in general not to optimal results regarding the Lebesgue constant. It turns out that for a quasi-optimal behavior of the Lebesgue constant, the entries of the Vandermonde matrix have to be normalized properly. If this normalization is adopted, our theoretical and numerical findings for the Lebesgue constant of the studied subfamilies of Fekete segments indeed show an asymptotic logarithmic growth. 

\section{Interpolation operators} \label{sect:interpolators}
We discuss and compare three kinds of polynomial interpolation operators $\Pi_r$. Based on $r$ data values, the operators $\Pi_r$ map functions defined on the reference interval $ I = [-1,1]$ to univariate polynomials $ \P_{r-1} $ of degree $ r-1 $. The choice of the reference interval is not restrictive and is thus set at convenience of computations. These interpolation operators are also projectors onto $ \P_{r-1} $, that is, they act as the identity operator when restricted to $ \P_{r-1} $. As a consequence of the Kharshiladze-Lozinski theorem, their operator norm with respect to the uniform norm is bounded from below (see \cite[Chap. 6, Sec. 5]{Cheney82} and \cite[Chap. 3, Thm. 1]{CheneyLight2000}) by 
\begin{equation} \label{eq:lowerboundLeb}
\Vert \Pi_r \Vert_{\op} \geq \frac{2}{\pi^2} \ln (r-1) - \frac{1}{2}.
\end{equation}
This has several relevant consequences for the interpolation with polynomials. It implies that there are continuous functions $f$ such that the interpolant $\Pi_r f$ does not converge towards $f$ for increasing $r$. It also implies that polynomial interpolation in general tends to have an increasing numerical conditioning as $r$ gets large. A logarithmic growth of the operator norm $\Vert \Pi_r \Vert_{\op}$ is the best that we can hope for such a polynomial interpolation problem.

\subsection{The nodal interpolator} \label{sect:nodalinterpolator}

Nodal polynomial interpolation is a celebrated tool of numerical analysis \cite{Cheney82} that uses function evaluations as input information for the interpolation. For a node set $\mathcal{X} = \{\xi_i\}_{i=1}^r$, it associates to any continuous function $ f \in C (I) $ a polynomial $ p_{r-1} \in \P_{r-1} $ of degree $ r-1 $ such that 
\begin{equation} \label{eq:nodalinterp}
	f(\xi_i) = p_{r-1}(\xi_i) \quad \text{ for } \;i \in \{1, \ldots, r\}.
\end{equation}
The interpolating polynomial $ p_{r-1} (x) = \Pi_{r}(\mathcal{X}) f(x) $ is unique when the nodes $\{ \xi_i \}_{i=1}^r $ are pairwise distinct; however, the proximity of $ p_{r-1} $ to $ f $, usually measured in terms of the uniform error $ \Vert f - p_{r-1} \Vert$ on $I$, depends strongly on the choice of the node set $\mathcal{X}$, and convergence with respect to $ r $ might not be granted, as for instance visible in the celebrated example given by Runge \cite{Runge}. 

In nodal polynomial interpolation one has the convenient representation
\begin{equation} \label{eq:ExplicitNodal}
\Pi_{r}(\mathcal{X}) f(x) = \sum_{i=1}^r f(\xi_i) \ell_{\xi_i} (x) ,
\end{equation}
where $ \{ \ell_{\xi_i} \}_{i=1}^r $ are the \emph{Lagrange basis polynomials} satisfying the duality relation $ \ell_{\xi_j} (\xi_i) = \delta_{i,j} $. They can be obtained from any convenient basis $ \{ B_j \}_{j=1}^r $ of $ \P_{r-1} $ by solving a linear system. 
In fact, defining the Vandermonde matrix $V({\mathcal{X}}) \in \R^{r \times r}$ in terms of the entries
\begin{equation} \label{eq:NodalVandermonde}
	V_{i,j} ({\mathcal{X}}) \doteq B_{j} (\xi_i), \quad i,j \in \{1, \ldots, r\},
\end{equation}
and the matrix $ V_k ({\mathcal{X}})(x) $ by replacing the $ k $-th row of $ V ({\mathcal{X}})$ by the vector $ (B_1(x), \ldots, B_r (x)) $, then, by Cramer's rule, one immediately gets
\begin{equation*} \label{eq:ellxi}
	\ell_{\xi_k} (x) = \frac{\det V_k ({\mathcal{X}})(x)}{\det V ({\mathcal{X}})}  .
\end{equation*}

\subsection{The segmental interpolator} \label{sect:segmentalinterpolator}

Segmental polynomial interpolation (or \emph{histopolation} \cite{Robidoux2006}) is a generalisation of nodal interpolation in which the point evaluations in \eqref{eq:nodalinterp} are replaced by integrals over segments $s_i = [\alpha_i,\beta_i] \subset I$, $i \in \{1, \ldots, r\}$, such that the respective histopolation condition reads as
\begin{equation}\label{eq:segmentalinterp}
	\int_{s_i} f(x) \de x  = \int_{s_i} p_{r-1}(x) \de x \quad \text{ for }\; i \in \{1, \ldots, r\} .
\end{equation}
This generalized interpolator emerges from physically-oriented problems \cite{Bossavitbook} and relaxes the regularity of the interpolating function, requiring $ f $ to be only essentially bounded on $ I $. For a set of segments $\mathcal{S} \doteq \{ s_i \}_{i=1}^r $ in the interval $I$, the interpolating (or \emph{histopolating}) polynomial $ p_{r-1}(x) = \Pi_{r}(\mathcal{S}) f(x) $ 
can be represented as
\begin{equation} \label{eq:segmentalexplicit}
	\Pi_r(\mathcal{S}) f(x) = \sum_{i=1}^{r} \left( \int_{s_i} f(x) \de x \right) \emm_{s_i} (x),
\end{equation}
where its cardinal basis functions satisfy the duality request
\begin{equation} \label{eq:dualityemm}
	\int_{s_i} \emm_{s_j}(x) \de x = \delta_{i,j}.
\end{equation}
Sufficient conditions for the set of segments $ \mathcal{S} = \{ s_i \}_{i=1}^r $ to give a unique polynomial interpolant $ \Pi_r(\mathcal{S}) f(x) $ are discussed in the literature, see \cite{Bojanov,BruniErb,Robidoux2006}. 
In particular, if the segments in $ \mathcal{S} $ are disjoint in measure, the mean value theorem can be applied to guarantee existence and uniqueness of the interpolation problem, see \cite[Proposition 3.1]{BruniErb} or \cite[Theorem 4.1]{Robidoux2006}. Under this hypothesis, the mean value theorem makes it also easy to localise the zeros of the Lagrange basis functions $ \{ \emm_{s_i} \}_{i=1}^{r} $. 

\begin{proposition} \label{prop:localisaationbases}
	Let $ \{ s_i \}_{i=1}^{r} $ be a unisolvent collection of disjoint segments and let $ \{ \emm_{s_i} \}_{i=1}^{r} $ be the corresponding Lagrange basis polynomials. Then
	\begin{itemize}
		\item[(i)] $ \emm_{s_i}(x) $ has $ 1 - \delta_{i,j} $ zeros in the segment $ s_j $, $j \in \{1, \ldots, r\}$.
	\end{itemize}
	As a consequence,
	\begin{itemize}
		\item[(ii)] $ \emm_{s_j} (x) > 0 $ for all $ x \in s_j $.
	\end{itemize}
\end{proposition}

{\setlength{\parindent}{0cm}
\textbf{Proof.}}
	By \eqref{eq:dualityemm} and the mean value theorem, $ \emm_{s_j}(x) $ has at least one zero in $ s_i $ for $ i \ne j $. It hence has $ r-1 $ zeros. Since it is a degree $ r-1 $ polynomial and it is not the zero polynomial, it has only these zeros. This proves $ (i) $. Since $ \emm_{s_j}(x) $ has no zeros in $ s_j $ and positive integral over $ s_j $, this also proves $ (ii) $.
\QED

In contrast to the nodal case, an explicit representation of the Lagrange functions $ \emm_{s_i} $ is known only in a few particular cases, for instance, when considering concatenated segments \cite{Gerritsma} or when all segments share the same left (or right) endpoint, see \cite{BruniErb}. Nevertheless, as in the nodal case a general recipe for a representing formula is obtained by fixing a convenient basis $ \{ B_j \}_{j=1}^r $ of $ \P_{r-1} $ and defining the Vandermonde matrix
\begin{equation} \label{eq:VandermondeSegments}
	\widehat{V}_{i,j} (\mathcal{S}) \doteq \int_{s_i} B_j (x) \de x .
\end{equation}
Letting $ \widehat{V}_k({\mathcal{S}})(x) $ denote the matrix obtained from $ \widehat{V} ({\mathcal{S}}) $ by replacing the $ k $-th row by the vector $ (B_1(x), \ldots, B_r (x)) $, again by Cramer's rule one gets
\begin{equation*} %\label{eq:ellxi}
	\emm_{s_k} (x) = \frac{\det \widehat{V}_k ({\mathcal{S}})(x)}{\det \widehat{V} ({\mathcal{S}})}  .
\end{equation*}

\subsection{A combined nodal-segmental interpolator based on function averages} \label{sect:averagedinterpolator}

If the lengths of the segments $ s_i $ are strictly positive, both sides of \eqref{eq:segmentalinterp} can be divided by the length $ |s_i| $ to obtain the equivalent interpolation condition
\begin{equation} \label{eq:averagedinterp}
	\frac{1}{|s_i|} \int_{s_i} f(x) \de x  = \frac{1}{|s_i|} \int_{s_i} p_{r-1}(x) \de x \quad \text{ for } \; i \in \{1, \ldots, r\}.
\end{equation}
The resulting interpolating polynomial $p_{r-1}$ does not change, however the Lagrange basis functions associated with \eqref{eq:averagedinterp} are now defined using the additional normalization
\begin{equation} \label{eq:dualityavg}
\frac{1}{|s_i|} \int_{s_i} \ell_{s_j}(x) \de x = \delta_{i,j}, 
\end{equation}
and, following the recipe offered in Section \ref{sect:segmentalinterpolator}, they can be retrieved as follows. Fix a convenient basis $ \{ B_j \}_{j=1}^r $ of $ \P_{r-1} $ and define the Vandermonde matrix
\begin{equation} \label{eq:VandermondeAveraged}
	V_{i,j} (\mathcal{S}) \doteq \frac{1}{|s_i|}\int_{s_i} B_j (x) \de x .
\end{equation}
Letting $ V_k ({\mathcal{S}})(x) $ denote the matrix obtained from $ V ({\mathcal{S}}) $ by replacing the $ k $-th row by the vector $ (B_1(x), \ldots, B_r (x)) $, Cramer's rule implies the formula
\begin{equation} \label{eq:ellsk}
	\ell_{s_k} (x) = \frac{\det V_k ({\mathcal{S}})(x)}{\det V ({\mathcal{S}})}  .
\end{equation}

Using the new cardinal functions \eqref{eq:ellsk}, the interpolator $ \Pi_{r}(\mathcal{S}) f(x) $ determined by the conditions \eqref{eq:averagedinterp} can be represented as
\begin{equation} \label{eq:explicitaveragedinterp}
	\Pi_r(\mathcal{S}) f(x) = \sum_{i=1}^{r} \left( \frac{1}{|s_i|}\int_{s_i} f(x) \de x \right) \ell_{s_i}(x) .
\end{equation}
Matching the representation \eqref{eq:explicitaveragedinterp} with the previously derived non-normalized description \eqref{eq:segmentalexplicit} of the interpolation operator, we obtain the relation 
\begin{equation} \label{eq:relationshipavgseg}
	\ell_{s_i}(x) = |s_i| \emm_{s_i}(x), \quad i \in \{1, \ldots, r\},
\end{equation}
between the Lagrange basis functions in the two descriptions. 

Under the hypothesis that the function $ f $ is continuous, the mean value theorem guarantees that the interpolation conditions \eqref{eq:averagedinterp} remain meaningful even if some or all of the segments $ s_i $ collapse to single points $\xi_i$. When allowing such limits, the conditions \eqref{eq:averagedinterp} give rise to a combined nodal-segmental interpolation problem based on function information on a set $\mathcal{A} = \{a_i\}_{i = 1}^r$, where the supports $a_i$ in $\mathcal{A}$ can either be single nodes or segments. We will distinguish between the following three cases. 
\begin{itemize}
\item[(A1)] \label{item:A1} \emph{Pure nodal interpolation}. All elements $ a_i \in \mathcal{A}$ are single points in $I$.
	\item[(A2)] \label{item:A2} \emph{Pure segmental interpolation}. All segment lengths $ |a_i| $ are strictly positive. If we define the diagonal matrix 
 $$ \mathbf{N} \doteq \begin{pmatrix}
     |a_1| & & \\ & \ddots & \\ && |a_r|
 \end{pmatrix}, $$ and compare the two Vandermonde matrices \eqref{eq:VandermondeSegments} and \eqref{eq:VandermondeAveraged}, we get the relation
 $ V (\mathcal{A}) = \mathbf{N}^{-1} \widehat{V}(\mathcal{A})$ for the two different normalizations of the input data. Further, the Lagrange functions for the two normalization possibilities satisfy $\ell_{a_i}(x) = |a_i| \emm_{a_i}(x)$, as in \eqref{eq:relationshipavgseg}.
\item[(A3)] \label{item:A3} \emph{Mixed nodal-segmental interpolation}. In this case, some of the elements $a_i \in \mathcal{A}$ satisfy $| a_i | = 0 $ and some $|a_i| > 0$. 
Separating the nodal and the segmental contributions, we may then expand the interpolation operator \eqref{eq:explicitaveragedinterp} as
\begin{equation*} \label{eq:hybridinterp}
	\Pi_r(\mathcal{A}) f(x) = \sum_{i = 1}^r \left( \frac{1}{|a_i|}\int_{a_i}\!\! f(x) \de x \right) \ell_{a_i}(x)  = \sum_{|a_i| = 0} \!\! f(a_i) \ell_{a_i}(x) + \sum_{|a_i| > 0} \! \left( \int_{a_i} \!\! f(x) \de x \right) \emm_{a_i}(x)  .
\end{equation*}
Here, for $|a_i| = 0$, the integral $\frac{1}{|a_i|}\int_{a_i}\!\! f(x) \de x$ denotes the point evaluation $f(a_i)$ of $f$ at the node $a_i$. For nodes $a_i$ we also have $ \ell_{a_i} (a_i) = 1 $ even if other constitutive conditions might be of the form $ \int_{a_j} \ell_{a_i}(x) \de x = 0 $.
\end{itemize}

\begin{remark} \label{prop:allshrinking}
We say that a sequence of segments $a_i^{(n)}$ converges to a segment $a_i$ for $n \to \infty$ if the endpoints $\alpha_i^{(n)}$ and $\beta_i^{(n)}$ of $a_i^{(n)}$ converge to the endpoints $\alpha_i$ and $\beta_i$ of $a_i$. As the Vandermonde determinants $V(\mathcal{A})$ and $V_k(\mathcal{A})(x)$ are continuous functions in terms of the parameters $\alpha_i$ and $\beta_i$ we get $\lim_{n\to \infty} V(\mathcal{A}^{(n)}) = V(\mathcal{A})$ and $\lim_{n\to \infty} V_k(\mathcal{A}^{(n)})(x) = V_k(\mathcal{A})(x)$ and therefore also 
\[ \lim_{n \to \infty} \ell_{a_i^{(n)}}(x) = \ell_{a_i}(x).\]
In particular, if in the limit $n \to \infty$ the left and right endpoints are indentical such that 
$\lim_{n \to \infty} \alpha_i^{(n)} = \lim_{n\to \infty} \beta_i^{(n)} = \xi_i$, 
we get as a limit of the sequence $a_i^{(n)}$ a single node $\xi_i$ and the nodal Lagrange functions $\ell_{\xi_i}(x)$ as a limit of the Lagrange functions with shrinking segments. This justifies the unified notation for the normalized Lagrange functions in the segmental and the nodal setting and shows that the nodal interpolation operator \eqref{eq:ExplicitNodal}  is a limiting case of the segmental interpolation operator \eqref{eq:explicitaveragedinterp}. For the non-normalized Lagrange functions $\emm_{a_i}(x)$ an additional normalization with the length of $a_i$ is required in order to be able to pass to the limiting case of segments that shrink to a single node, see \eqref{eq:relationshipavgseg}.
\end{remark}

As the nodal interpolation operator \eqref{eq:ExplicitNodal} or the segmental interpolation operator \eqref{eq:segmentalexplicit} also the combined nodal-segmental interpolator is a projector. 

\begin{proposition} \label{prop:avgprojector}
For a unisolvent set of supports $\mathcal{A}$ consisting of $r$ nodes and segments, the nodal-segmental interpolator satisfies $ \Pi_{r}(\mathcal{A})^2 f(x) = \Pi_{r}(\mathcal{A}) f(x) $.
\end{proposition}

{\setlength{\parindent}{0cm}
\textbf{Proof.}}
	From \eqref{eq:dualityavg} %the splitting in \eqref{eq:hybridinterp} we can 
    we deduce that $ \Pi_r(\mathcal{A}) \ell_{a_i} = \ell_{a_i} $. 
	The linearity of $ \Pi_{\mathcal{A}} $ then yields
	\begin{align*}
		\Pi_r(\mathcal{A}) \left( \Pi_r(\mathcal{A}) f(x) \right) & = \Pi_r(\mathcal{A}) \left( \sum_{i = 1}^r \left( \frac{1}{|a_i|}\int_{a_i} f(x) \de x \right) \ell_{a_i}(x) \right) \\
        & = \sum_{i = 1}^r \left( \frac{1}{|a_i|}\int_{a_i} f(x) \de x \right) \Pi_r(\mathcal{A}) \ell_{s_i}(x)  \\
		& = \sum_{i = 1}^r \left( \frac{1}{|a_i|}\int_{a_i} f (x) \de x \right) \ell_{a_i}(x) = \Pi_r(\mathcal{A}) f(x),
	\end{align*}
    where, as before, the average value $\frac{1}{|a_i|}\int_{a_i} f (x) \de x$ corresponds to the function evaluation $f(a_i)$ in the limit case $|a_i| = 0$. 
	This concludes the proof.
\QED

We can also merge nodal and segmental results for unisolvence of polynomial interpolation to get sufficient conditions for the well-posedness of the combined interpolator $ \Pi_r(\mathcal{A})$.

\begin{definition} \label{def:unisolvence}
A set of supports $ \mathcal{A} $ is said to be regular if it splits as $ \mathcal{A} = \mathcal{X} \cup \mathcal{S} $ and
	\begin{itemize}
		\item[(i)] $ \mathcal{X} $ contains pairwise distinct nodes;
		\item[(ii)] $ \mathcal{S} $ contains segments that overlap at most in their endpoints;
		\item[(iii)] points in $ \mathcal{X} $ are not contained in the interior of any element of $ \mathcal{S} $.
	\end{itemize}
\end{definition}

Regular sets are unisolvent: by the mean value theorem, one sees that a polynomial of degree $r-1$ whose averages vanish on a regular set has  $ r $ distinct zeros, so it is the zero polynomial. Hence, a regular set $ \mathcal{A} $ provides a well-posed interpolator $ \Pi_r(\mathcal{A}) $ in the space of polynomials $\P_{r-1}$, see \cite[Theorem 6.5]{Robidoux2006}.

\section{Lebesgue constants} \label{sect:lebconst}
Lebesgue's lemma states that, whenever $ \Pi_r : V \to \P_{r-1} \subseteq V $ is a linear projector from a normed space $V$ onto the space of polynomials of degree $r-1$, we have the inequality (see \cite[Chap. 2]{Davis75})
$$ \Vert f - \Pi_r f \Vert_V \leq (1 + \Vert \Pi_r \Vert_{\op}) \inf_{p \in \P_{r-1}} \Vert p - \Pi_r f \Vert_V $$
provided that $ \P_{r-1} $ inherits the norm from $ V $. The induced operator norm $ \Vert \Pi_r \Vert_{\op} $ therefore gives a measure for the numerical conditioning of the interpolation problem and for the general quality of the interpolator. By choosing an appropriate norm, in our case the uniform norm, the quantity $ \Vert \Pi_r \Vert_{\op} $ can be characterised in terms of the involved nodes and segments, and thus helps in selecting a reliable family of supports \cite{AnaFra}.

\subsection{The nodal Lebesgue constant}

The Lebesgue constant associated with a collection of nodes $ \mathcal{X} = \{ \xi_i\}_{i=1}^r $ is given as
\begin{equation} \label{eq:nodalLeb}
\Lambda_r (\mathcal{X}) \doteq \sup_{x \in I} \sum_{i=1}^{r} | \ell_{\xi_i} (x) |. 
\end{equation}
When $ C(I) $ is endowed with the sup-norm and $ \Pi_r(\mathcal{X}): C (I) \to \P_{r-1} $ is the nodal interpolator introduced in Section \ref{sect:nodalinterpolator}, one has the identity
\begin{equation} \label{eq:nodalnorm}
\Vert \Pi_r(\mathcal{X}) \Vert_{\op} = \Lambda_r (\mathcal{X}),
\end{equation} 
i.e., the operator norm $\Vert \Pi_r(\mathcal{X}) \Vert_{\op}$ corresponds to $\Lambda_r (\mathcal{X})$. 
The Lebesgue constant $\Lambda_r (\mathcal{X})$ is highly sensitive to the choice of the set $ \mathcal{X} $, and asymptotic growth behaviours ranging from logarithmic \cite{Sundermann} to exponential growth \cite{Trefethen} have been observed. In view of the lower bound \eqref{eq:lowerboundLeb}, the identity \eqref{eq:nodalnorm} gives rise to the search of  the node sets $\mathcal{X}$ for which the constants $\Lambda_r (\mathcal{X})$ grows logarithmically in $r$. For a survey on relevant results we refer to \cite{Ibrahimoglu}.

\subsection{The segmental Lebesgue constant}

The Lebesgue constant associated with a collection of segments $ \mathcal{S} = \{ s_i \}_{i = 1}^r $ is given as
\begin{equation} \label{eq:segmentalLeb} \Lambda_{r} (\mathcal{S}) \doteq \sup_{s \subseteq I} \frac{1}{|s|} \sum_{i=1}^{r} |s_i| \left\vert \int_s \emm_{s_i} (x) \de x \right\vert = \sup_{x \in I} \sum_{i=1}^{r} |s_i| | \emm_{s_i} (x) | .
\end{equation}
It was shown in \cite[Theorem $4.1$]{BruniErb} that, if the segments $ s_i $ overlap at most in their endpoints, the operator norm of the interpolator \eqref{eq:segmentalexplicit} satisfies
\begin{equation} \label{eq:segmentalnorm}
	\Vert \Pi_r(\mathcal{S}) \Vert_{\op} = \Lambda_r (\mathcal{S}) .
\end{equation}
This holds true either in the case $ \Pi_r (\mathcal{S}): C (I) \to \P_{r-1}$ or in the case $ \Pi_r(\mathcal{S}): L_\infty (I) \to \P_{r-1}$, with both spaces endowed with the sup-norm. In the more general case of overlapping segments, the quantity \eqref{eq:segmentalLeb} still offers an upper bound for $ \Vert \Pi_r (\mathcal{S}) \Vert_{\op} $.

Similarly to the nodal Lebesgue constant \eqref{eq:nodalLeb}, also the segmental Lebesgue constant \eqref{eq:segmentalLeb} has proved to be highly sensitive to the choice of the segments $ \mathcal{S} $. This was first observed numerically \cite{BruniEdges} and then proved in the following specific scenarios \cite{BruniErb}:
\begin{itemize}
	\item[(C1)] \emph{Concatenated segments}. Given a set of nodes $ -1 = \xi_0 < \ldots < \xi_{r} = 1 $, the concatenated segments are defined as $ s_i \doteq [\xi_{i-1}, \xi_{i}] $. Note that $ r+1 $ nodes are required.
	\item[(C2)] \emph{Segments with uniform arc-length}. Fixed an arc-radius $ 0 < \rho < \pi $ and a set of values $ 0<\tau_1 < \ldots < \tau_{r} < \pi$, we define the segments $ s_i \doteq [\cos(\tau_i + \rho), \cos(\tau_i - \rho)] $. Only $ r $ nodes $\tau_i$ and a parameter $\rho$ are required for the definition.
\end{itemize}
In both cases (C1) and (C2) a connection between the point set and the segmental set was observed. In particular, it was proved that the segmental Lebesgue constant of supports in the class (C1) with $ \xi_{i} = -1 + \frac{2i}{r} $ shows an exponential growth \cite[Eq. (5.3)]{BruniErb}, whereas that of segments in the class (C2) with $ \tau_i = \frac{2 (i-1) -1}{2r} \pi $ shows a logarithmic growth \cite[Corollary 5.7]{BruniErb}. In view of the bound \eqref{eq:lowerboundLeb} this latter case is optimal up to a constant.

\subsection{The Lebesgue constant for mixed nodal-segmental interpolation} \label{sect:averagedLeb}

When nodal evaluations and integrals of a function are combined as input data of the interpolator, we also want to have a unified description of the Lebesgue constants that contains the nodal \eqref{eq:nodalLeb} as well as the segmental form \eqref{eq:segmentalLeb}. To describe the norm of the combined nodal-segmental interpolation operator \eqref{eq:explicitaveragedinterp}, we therefore define
\begin{equation} \label{eq:averagedLeb}
	\Lambda_r(\mathcal{A}) \doteq \sup_{s \subset I} \frac{1}{|s|} \sum_{i=1}^{r} \left\vert \int_s \ell_{a_i}(x) \de x  \right\vert = \sup_{x \in I} \sum_{i=1}^{r} \left\vert \ell_{a_i}(x) \right\vert.
\end{equation}

It is obvious that this definition contains automatically the Lebesgue constant \eqref{eq:nodalLeb} in the nodal scenario (A1). By the relation \eqref{eq:relationshipavgseg} between the Lagrange functions $\ell_{a_i}$ and $\emm_{a_i}$ it is also easily seen that this definition is consistent with \eqref{eq:segmentalLeb} in the segmental scenario (A2).  

For a regular set $\mathcal{A}$ introduced in Definition \ref{def:unisolvence}, the operator norm of $ \Pi_r(\mathcal{A}) $ turns out to be identical to the Lebesgue constant \eqref{eq:averagedLeb}.

\begin{theorem} \label{thm:equaivalenceLebesgue}
	Let $ \P_{r-1} \subset C(I)$ be endowed with the sup-norm. If the set $ \mathcal{A}$ of $r$ supports (nodes and/or segments) is unisolvent for $\P_{r-1}$, we have the inequality $$ \Vert \Pi_r(\mathcal{A}) \Vert_{\op} \leq \Lambda_r (\mathcal{A}). $$ 
    If $ \mathcal{A} = \mathcal{S} \cup \mathcal{X} $ is regular (i.e. $\mathcal{A}$ satisfies Definition \ref{def:unisolvence}) we have equality
	\begin{equation*}
		\Vert \Pi_r(\mathcal{A}) \Vert_{\op} = \Lambda_r (\mathcal{A}) .
	\end{equation*}

\end{theorem}
%
%Before turning to the proof of Theorem \ref{thm:equaivalenceLebesgue}, we evidence the criticality of having a node in $ \mathcal{X} $ which is also an endpoint of a segments in $ \mathcal{S}$. This occurs:
%\begin{itemize}
%    \item[$(i)$] if an endpoint of an isolated edge in $ \mathcal{S} $ belongs to $ \mathcal{X} $;
%    \item[$(ii)$] if two edges in $ \mathcal{S} $ share a vertex, which is a node in $ \mathcal{X} $;
%    \item[$(iii)$] if two nodes in $ \mathcal{X} $ are the endpoints of an edge in $ \mathcal{S} $.
%\end{itemize}
%All the above situations force to introduce a parameter $ \varepsilon > 0 $, which is appropriately sent to zero to obtain the result. We propose a proof for the case $ (ii) $, which immediately extends to $ (i) $ and to $ (iii) $ with small effort.

{\setlength{\parindent}{0cm}
\textbf{Proof.}}
    We show that $ \Vert \Pi_r (\mathcal{A}) \Vert_{\op} \leq \Lambda_r (\mathcal{A}) $. Inserting the respective definitions and applying the triangular inequality, we get
    \begin{align*} \Vert \Pi_r (\mathcal{A}) \Vert_\op & = \sup_{\Vert f \Vert = 1} \Vert \Pi_r f \Vert = \sup_{\Vert f \Vert = 1} \sup_{x \in I } \vert \Pi_r f \vert = \sup_{\Vert f \Vert = 1} \sup_{x \in I } \left\vert \sum_{i=1}^r \left( \frac{1}{|a_i|} \int_{a_i} f \right) \ell_{a_i} (x) \right\vert \\
    & \leq \sup_{\Vert f \Vert = 1} \sup_{x \in I }  \sum_{i=1}^r \left\vert \frac{1}{|a_i|}  \int_{a_i} f  \right \vert\left\vert \ell_{a_i} (x) \right\vert \leq \sum_{i=1}^r \left\vert \ell_{a_i} (x) \right\vert = \Lambda_r (\mathcal{A} ),
    \end{align*}
    where the last inequality is granted by the mean value theorem, since $ \Vert f \Vert = 1 $ ensures that $ \frac{1}{|a_i|} \left\vert\int_{a_i} f \right\vert \leq 1 $ for each $ a_i \subset I $. This proves the first part of the statement.

    To prove that $ \Vert \Pi_r (\mathcal{A}) \Vert_\op \geq \Lambda_r $, we construct a continuous function $ f_\varepsilon \in C(I)$ such that $ \Vert \Pi_r f_\varepsilon \Vert \geq \Lambda_r (\mathcal{A}) - \varepsilon $ for any small $ \varepsilon > 0 $. With the notation of Definition \ref{def:unisolvence}, we split $ \mathcal{A} = \mathcal{X} \cup \mathcal{S}$ and recall that $ |a_i| = 0 $ if $ a_i \in \mathcal{X}$ and $ |a_i| > 0 $ if $ a_i \in \mathcal{S}$. We write $ a_i = [\alpha_i, \beta_i] $ with the convention that $ \alpha_i = \beta_i $ if $ |a_i|=0 $. To any $ a_i \in \mathcal{A} $ we associate the neighbourhoods
    $$ U_i \doteq \begin{cases}
        \{ a_i \}
        \quad \text{ if } \; a_i \in \mathcal{X}, \\
        \left(\alpha_i + \varepsilon, \beta_i - \varepsilon \right) \quad \text{ if } \; a_i \in \mathcal{S} ,
    \end{cases}
    \text{and} \quad  V_i \doteq \begin{cases}
        (\alpha_i - \varepsilon \eta /4, \beta_i + \varepsilon \eta/4)
        \quad \text{ if } \; a_i \in \mathcal{X}, \\
        \left(\alpha_i + \varepsilon/2, \beta_i - \varepsilon/2 \right) \quad \text{ if } \; a_i \in \mathcal{S} ,
    \end{cases}
    $$
    where $\eta = \min \{2,|a_i| > 0\}$. Then, $0 < \eta \leq 2$, and by the regularity of $ \mathcal{A} $, we have $ U_i \subset V_i $ for each $ i = 1, \ldots, r $ and $ V_i \cap V_j = \emptyset $ for $ i \ne j $.

    We now consider a collection of $ r $ piecewise linear functions $ \{ \varphi_i \}_{i=1}^r : I \to [0,1] $ such that
    $$ \varphi_i (x) = \begin{cases}
        1 \quad \text{ on } \quad U_i, \\
        0 \quad \text{ outside } \quad V_i .
    \end{cases} $$
    We have the inclusions $ a_i \subseteq U_i \subseteq V_i $ if $ a_i \in \mathcal{X} $ and $ U_i \subseteq V_i \subseteq a_i $ if $ a_i \in \mathcal{S} $. Hence, if $ | a_i | = 0 $, the function $ \varphi_i (x) $ is the only non-vanishing on $ a_i $ and 
    $ \frac{1}{|a_i|}\int_{a_i} \varphi_i (x) \de x = \varphi_i (a_i) = 1 \geq (1 - \varepsilon) $. If $ |a_i| > 0 $, up to three consecutive functions $ \varphi_{i-1}$, $ \varphi_i$ and $ \varphi_{i+1}$ may not vanish on $ a_i $. We have
\begin{equation} \label{eq:firsttechnical}
\int_{a_i} \varphi_i (x) \de x = \int_{V_i} \varphi_i (x) \de x \geq \int_{U_i} \varphi_i (x) \de x = (1-\varepsilon) |a_i| \end{equation}
    and
\begin{equation} \label{eq:secondtechnical}
\int_{a_i} \varphi_{i+1} (x) \de x = \int_{a_i \cap V_{i+1}} \!\!\! \varphi_{i+1} (x) \de x \leq \frac{\varepsilon |a_i|}{4}, \quad \int_{a_i} \varphi_{i-1} (x) \de x = \int_{a_i \cap V_{i-1}} \!\!\! \varphi_{i-1} (x) \de x \leq \frac{\varepsilon |a_i|}{4} 
\end{equation}
Since $ I $ is compact, we may let $ \bar{x} \doteq \arg \max_{x \in I} \sum_{i=1}^r \left\vert \ell_{a_i} (x) \right\vert $ and define
    $$ f_\varepsilon(x) \doteq \sum_{i=1}^r \mathrm{sgn} \left(\ell_{a_i} (\bar{x}) \right) \varphi_i (x) .$$
    Since $ V_i \cap V_j = \emptyset$ for $ i \ne j $, $ f_\varepsilon $ is continuous and $ \Vert f_\varepsilon \Vert = 1 $. We thus get
    {\allowdisplaybreaks \begin{align*}
        \Vert \Pi_r (\mathcal{A}) \Vert_\op & = \sup_{\Vert f \Vert = 1} \Vert \Pi_r f \Vert \geq \Vert \Pi_r f_\varepsilon \Vert  = \sup_{x \in I} \left\vert \sum_{i=1}^r \left( \frac{1}{|a_i|} \int_{a_i} f_\varepsilon (x) \de x \right) \ell_{a_i} (x) \right\vert \\
        & = \sup_{x \in I} \left\vert \sum_{i=1}^r \left( \frac{1}{|a_i|} \int_{a_i} \sum_{j=1}^r \mathrm{sgn} \left(\ell_{a_j} (\bar{x}) \right) \varphi_j (x) \de x \right) \ell_{a_i} (x) \right\vert \\
        & = \sup_{x \in I} \left\vert \sum_{i=1}^r \left( \frac{1}{|a_i|} \int_{a_i} \sum_{j=i-1}^{i+1} \mathrm{sgn} \left(\ell_{a_j} (\bar{x}) \right) \varphi_j (x) \de x \right) \ell_{a_i} (x) \right\vert .
    \end{align*}}
    We apply the triangular inequality to split the three summands:
    \begin{align*}
        \Vert \Pi_r (\mathcal{A})\Vert_{\op} & = \sup_{x \in I} \left\vert \sum_{i=1}^r \left( \frac{1}{|a_i|} \int_{a_i} \sum_{j=i-1}^{i+1} \mathrm{sgn} \left(\ell_{a_j} (\bar{x}) \right) \varphi_j (x) \de x \right) \ell_{a_i} (x) \right\vert \\
        & \geq \sup_{x \in I} \left\vert \sum_{i=1}^r \mathrm{sgn} \left(\ell_{a_i} (\bar{x}) \right) \left( \frac{1}{|a_i|} \int_{a_i}  \varphi_i (x) \de x \right) \ell_{a_i} (x) \right\vert \\ & - \sup_{x \in I} \left\vert \sum_{i=1}^r \mathrm{sgn} \left(\ell_{a_i} (\bar{x}) \right) \left( \frac{1}{|a_i|} \int_{a_i}  \varphi_{i-1} (x) \de x \right) \ell_{a_i} (x) \right\vert \\ & - \sup_{x \in I} \left\vert \sum_{i=1}^r \mathrm{sgn} \left(\ell_{a_i} (\bar{x}) \right) \left( \frac{1}{|a_i|} \int_{a_i}  \varphi_{i+1} (x) \de x \right) \ell_{a_i} (x) \right\vert .
        \end{align*}
        By \eqref{eq:firsttechnical}, we may bound the first term as
        \begin{align*} 
        & \sup_{x \in I} \left\vert \sum_{i=1}^r \mathrm{sgn} \left(\ell_{a_i} (\bar{x}) \right) \left( \frac{1}{|a_i|} \int_{a_i}  \varphi_i (x) \de x \right) \ell_{a_i} (x) \right\vert \\ & \;\geq (1-\varepsilon) \sup_{x \in I} \left\vert \sum_{i=1}^r \mathrm{sgn} \left(\ell_{a_i} (\bar{x}) \right) \ell_{a_i} (x) \right\vert
        = (1-\varepsilon) \sup_{x \in I} \sum_{i=1}^r \left\vert \ell_{a_i} (x) \right\vert = (1-\varepsilon) \Lambda_r (\mathcal{A}) .
        \end{align*}
        By \eqref{eq:secondtechnical} and the triangle inequality, we may also bound the remaining two quantities:%. We show it just for $ \varphi_{i-1} (x) $, the other is identical:
        \begin{align*}
        & \sup_{x \in I} \left\vert \sum_{i=1}^r \mathrm{sgn} \left(\ell_{a_i} (\bar{x}) \right) \left( \frac{1}{|a_i|} \int_{a_i}  \varphi_{i-1} (x) \de x \right) \ell_{a_i} (x) \right\vert \\
        & \; \leq \sup_{x \in I}  \sum_{i=1}^r \left\vert \left( \frac{1}{|a_i|} \int_{a_i}  \varphi_{i-1} (x) \de x \right) \ell_{a_i} (x) \right\vert 
        \leq \sup_{x \in I} \frac{\varepsilon}{4} \sum_{i=1}^r \left\vert \ell_{a_i} (x) \right\vert = \frac{\varepsilon}{4} \Lambda_r (\mathcal{A}) .
        \end{align*}
        The bound for $ \varphi_{i+1} (x) $ is identical. Gathering all the above estimates, we obtain
        $$\Vert \Pi_r (\mathcal{A}) \Vert_\op \geq  (1-\varepsilon) \Lambda_r (\mathcal{A}) - 2 \frac{\varepsilon}{4} \Lambda_r (\mathcal{A})
        = \left(1 -\frac{3}{2}\varepsilon \right) \Lambda_r (\mathcal{A}) .$$
    Since $ \Lambda_r (\mathcal{A}) $ is finite, the claim follows by defining $ \varepsilon' = \frac{3}{2}\varepsilon \Lambda_r (\mathcal{A}) $ and letting $ \varepsilon' \to 0 $.
\QED
 
\begin{remark}
If the regular set $ \mathcal{A} $ is in the class (A1), i.e. $  \mathcal{A} = \mathcal{X} $, Theorem \ref{thm:equaivalenceLebesgue} implies \eqref{eq:nodalnorm}. If the regular set $ \mathcal{A}$ is in the class (A2), i.e. $  \mathcal{A} = \mathcal{S} $, Theorem \ref{thm:equaivalenceLebesgue} implies \eqref{eq:segmentalnorm}. \end{remark}
 
\section{Fekete problems} \label{sect:feketeprobs}

The relation between the Lebesgue constant \eqref{eq:averagedLeb} and the operator norm of the interpolation operator shown in Theorem \ref{thm:equaivalenceLebesgue} suggests that a valuable set $\mathcal{A}$ for polynomial interpolation contains supports $a_i$ that keep the absolute values of the Lagrange functions $\ell_{a_i}(x)$ low. One strategy for this consists in identifying sets $\mathcal{A}^{\mathrm{Fek}}$ that maximise the absolute value of the determinant of the Vandermonde matrices $|\det V(\mathcal{A})|$ or $|\det \widehat{V}(\mathcal{A})|$. 
These sets, which are easily proven to be blind to the selected polynomial basis, will be referred to as \emph{Fekete sets}. For nodal interpolation on the interval, they are well-known in the literature \cite{Fejer1932}. 

\subsection{The Lagrange basis functions for Fekete nodes and segments}
We prove first of all that if one chooses a set $ \mathcal{A}^{\mathrm{Fek}} $ of supports  that maximise the absolute value of the Vandermonde determinant $ |\det V ({\mathcal{A}})| $, the corresponding Lagrange basis functions are small with respect to the sup-norm.

\begin{lemma} \label{lem:boundzeronormaveraged}
	Assume that the set of supports $ \mathcal{A}^{\mathrm{Fek}} = \{a_1, \ldots, a_r\} $ satisfies $ |\det V({\mathcal{A}})| \leq |\det V ({\mathcal{A}^{\mathrm{Fek}}})| $ for all sets $ \mathcal{A} $ of size $r$. Then $$ \Vert \ell_{a_i} \Vert = \sup_{x \in I} |\ell_{a_i}(x) | = \sup_{a \subseteq I} \frac{1}{|a|}\left \vert \int_a \ell_{a_i} (x) \de x \right \vert = 1 $$
	for all $ i \in \{ 1, \ldots, r\}$.
\end{lemma}

{\setlength{\parindent}{0cm}
\textbf{Proof.}}
	Suppose the segments $ \mathcal{A}^{\mathrm{Fek}} = \left\{ a_i \right\}_{i=1}^{r} $ maximise the determinant of $ V (\mathcal{A}) $ in the above sense. 
	Expanding the representation \eqref{eq:ellsk} via the Laplace theorem along the $i$-th row, integrating and dividing by $ |a| $, and using the properties of the determinant, one obtains
	$$ \frac{1}{|a|}\int_a \det V_i (\mathcal{A}^{\mathrm{Fek}})(x) \de x =  \det \begin{pmatrix}
		\frac{1}{|a_1|}\int_{a_1} B_1(x) \de x & \ldots & \frac{1}{|a_1|}\int_{a_1} B_{r}(x) \de x \\
		\vdots & \vdots & \vdots \\
		\frac{1}{|a|}\int_a B_1(x) \de x & \ldots & \frac{1}{|a|} \int_a B_{r} (x) \de x \\
		\vdots & \vdots & \vdots \\
		\frac{1}{|a_r|}\int_{a_{r}} B_1(x) \de x & \ldots & \frac{1}{|a_r|}\int_{a_{r}} B_{r}(x) \de x 
	\end{pmatrix} ,$$
	where the average $ \frac{1}{|a|} \int_a B_j(x) \de x $ is understood as a point evaluation $B_j(a)$ if $ |a| = 0$.
	
	We suppose, by contradiction, that $ \frac{1}{|a|} \left| \int_a \det V_i (\mathcal{A}^{\mathrm{Fek}})(x) \de x \right|  > |\det V (\mathcal{A}^{\mathrm{Fek}})| $ for some $ a \not\in \mathcal{A} $. Then, the set $ \mathcal{\tilde{A}} = \left\{ a_j \right\}_{j\ne i} \cup a $ satisfies $|\det V(\mathcal{\tilde{A}})| > |\det V(\mathcal{A}^{\mathrm{Fek}})| $, which is in contrast to the assumption that $\mathcal{A}^{\mathrm{Fek}}$ maximizes the  absolute value of the Vandermonde determinant. As a consequence, we have $$ \frac{1}{|a|} \left\vert \int_a \det V_i (\mathcal{A}^{\mathrm{Fek}}) (x) \de x \right\vert \leq \left\vert \det V (\mathcal{A}^{\mathrm{Fek}}) \right\vert, $$ 
    and hence
	$$ 1 \geq \frac{1}{|a|} \frac{\left\vert \int_a \det V_i (\mathcal{A}^{\mathrm{Fek}})(x) \de x \right\vert}{|\det V(\mathcal{A}^{\mathrm{Fek}})|} =  \frac{1}{|a|} \left\vert  \int_a \ell_{a_i}(x) \de x \right\vert$$
 for all possible points or segments $a$ in $I$. Here, equality is attained for $a = a_i$. 
\QED

When restricting the sets $\mathcal{A}$ in Lemma \ref{lem:boundzeronormaveraged} to the particular nodal scenario (A1), we get an analog and well-known result for the nodal Fekete points \cite{Fejer1932b}.

\begin{lemma} \label{lem:normboundnodal}
	Let $ \mathcal{X}^{\mathrm{Fek}} $ be a set of nodes such that $ |\det V ({\mathcal{X}})| \leq |\det V ({\mathcal{X}^{\mathrm{Fek}}})| $ for each node set $ \mathcal{X}$ with $r$ nodes. Then
	$$ \Vert  \ell_{\xi_i} \Vert = \max_{x \in I} \left\vert \ell_{\xi_i} (x) \right\vert = 1 $$
	for each $ i \in \{1, \ldots, r\}$.
\end{lemma}

If instead of considering sets that maximize $ |\det V ({\mathcal{A}})|$ one looks for Fekete segments $\widehat{\mathcal{S}}^{\mathrm{Fek}}$ that maximise the non-normalized determinant $ |\det \widehat{V} ({\mathcal{S}})| $, then one can prove that the corresponding basis functions $\emm_{s_i}$ have small integrals. The respective proof follows the lines of the proof of the more general Lemma \ref{lem:boundzeronormaveraged}.

\begin{lemma} \label{lem:integralsegmentalLagrange}
	Let $ \widehat{\mathcal{S}}^{\mathrm{Fek}} $ be such that $ |\det \widehat{V}({\mathcal{S}})| \leq |\det \widehat{V} ({\widehat{\mathcal{S}}^{\mathrm{Fek}}}) $ for all sets $ \mathcal{S} $ of size $r$.
	Then, for each segment $ s \subseteq I $, we have
	$$ \left \vert \int_s \emm_{s_i} (x) \de x \right \vert \leq 1 $$
	for all $ i \in \{ 1, \ldots, r\} $.
\end{lemma}

\subsection{The nodal Fekete problem}

We denote by $ \mathcal{X}^{\mathrm{Fek}} = \{\xi_1, \ldots, \xi_r \} $ the collection of nodes that maximises $ |\det V (\mathcal{X})| $. If such a collection exists, it does not depend on the chosen basis for $ \P_{r-1} $, see \cite{Bos1D}, and its points are typically called \emph{Fekete nodes}. Lemma \ref{lem:normboundnodal} immediately yields the following result.
\begin{lemma} \label{lem:estimateFeketeNodes}
	One has
	$$ \Lambda_r (\mathcal{X^{\mathrm{Fek}}}) \leq r.$$
\end{lemma}
By using the monomial basis $ \mathfrak{M} = \{1, x, \ldots, x^{r-1}\}$ in the Vandermonde matrix \eqref{eq:NodalVandermonde}, one obtains the explicit formula for the determinant
\begin{equation} \label{eq:detVDMproduct}
	\det V^{\mathfrak{M}} (\mathcal{X}) = \det \begin{pmatrix}
		(\xi_1)^0 & \ldots & (\xi_1)^{r-1} \\
		\vdots & \ddots & \vdots \\
		(\xi_{r})^0 & \ldots & (\xi_{r})^{r-1}
	\end{pmatrix} = \prod_{1 \leq i<j \leq r} \left( \xi_j - \xi_i \right) ,
\end{equation}
and hence the nodal Fekete problem may be stated in the following neat formulation.  
\begin{problem} \label{prob:Fek}
	Find a collection of nodes $ -1 \leq \xi_1 < \ldots < \xi_{r} \leq 1 $ that maximise the quantity
	\begin{equation*} \label{eq:productFekete}
		\prod_{1 \leq i<j \leq r} \left\vert \xi_j - \xi_i \right\vert .
	\end{equation*}
\end{problem}

For the uniqueness of the solution of Problem \ref{prob:Fek} the fact that all $\xi_i$ are contained in the inteval $I$ is central. A further assumption to obtain a unique maximizer of Problem \ref{prob:Fek} is the point symmetry of the nodes with respect to the center of $I$. It is well-known that points that solve the Fekete problem on the interval $I$ coincide with the \emph{Legendre-Gauss-Lobatto} nodes \cite{Fejer1932}. For the $d$-dimensional cube cube $I^{d}$ it is further known that the respective Fekete problem is solved by the tensor-product Legendre-Gauss-Lobatto points \cite{Bos2001}.

\begin{proposition} \label{lem:legendrezeros}
 Let $ P_{r-1} (x) $ be the Legendre polynomial of degree $r-1$. The roots $\xi_i^{(LL)}$, $i \in \{1, \ldots, r\}$ of the polynomial $ (1-x^2) P_{r-1}' (x) $ (known as Legendre-Gauss-Lobatto nodes) solve the nodal Fekete problem \ref{prob:Fek} in the interval $ I =  [-1, 1] $.   
\end{proposition} 

As shown in \cite{Fejer1932}, Legendre-Gauss-Lobatto points induce Lagrange basis functions that satisfy the inequality $$ \sum_{i=1}^{r} \vert \ell_{\xi_i} (x) \vert^2 \leq 1. $$ This is sometimes also referred to as \emph{Fej\'er condition}. This inequality in combination with the Cauchy-Schwarz inequality immediately implies the improved estimate 
\begin{equation} \label{eq:sqrtnodLeb}
	\Lambda_r (\mathcal{X^{\mathrm{Fek}}}) =  \sup_{x \in I} \sum_{i=1}^{r} \vert \ell_{\xi_i} (x) \vert \leq \sqrt{r} \sqrt{ \sum_{i=1}^{r} \vert \ell_{\xi_i} (x) \vert^2 } \leq \sqrt{r}
\end{equation}
for the Lebesgue constant of the Fekete nodes \cite{Fejer1932b}. This is, however, not the end of the story. With localisation techniques one shows that the Lebesgue constant for such points offers a logarithmic growth \cite{Sundermann} with unknown constant factors. These constants have been estimated numerically, for instance, in \cite{Hesthaven}, where it is conjectured that
\begin{equation*} \label{eq:LebLGL}
	\Lambda_r (\mathcal{X^{\mathrm{Fek}}}) \leq \frac{2}{\pi}\log(r) + c,
\end{equation*}
with $ c \approx 0.685 $. 
This, in view of \eqref{eq:lowerboundLeb}, ensures that Fekete nodes are quasi-optimal in terms of the growth of the Lebesgue constant. For a valuable survey on Fekete nodes, see \cite{BosSurvey}.

\subsection{The non-normalized Fekete problem for concatenated segments}

A relevant obstacle in the analytic resolution of the segmental Fekete problem is the larger number of parameters at play. In this section, we confine ourselves to the case (C1) of concatenated segments, where the number of free parameters gets reduced to $r+1$.  Proposition \ref{prop:extremashrinking} will treat this case afterwards in the normalized segmental framework.

\begin{problem} \label{prob:concatenated}
	Find concatenated segments $ \widehat{\mathcal{S}}^{\mathrm{Fek}} $ that maximise the determinant of the non-normalized Vandermonde matrix $\widehat{V}(\mathcal{S})$ given in \eqref{eq:VandermondeSegments}.
\end{problem}

The solution to Problem \ref{prob:concatenated} is intimately related to the nodal Fekete Problem \ref{prob:Fek}. Since we are assuming $ \widehat{\mathcal{S}}^{\mathrm{Fek}} $ to be concatenated, we may write the segments $s_i$ in $\widehat{\mathcal{S}}^{\mathrm{Fek}}$ in terms of their endpoints $\mathcal{X} = \{\xi_1, \ldots, \xi_{r+1}\}$ as $ s_i = [\xi_{i}, \xi_{i+1}] $, $i \in \{1, \ldots, r\}$. Choosing the monomial basis $ \mathfrak{M} = \{1, x, \ldots, x^{r-1}\}$ as representing system for the polynomials in $\P_{r-1}$, we can express the Vandermonde matrix of the segmental interpolation problem as
	$$ \widehat{V}^{\mathfrak{M}}(\mathcal{S}) = \begin{pmatrix}
	\xi_2-\xi_1 & \ldots & \frac{\xi_2^{r}-\xi_1^{r}}{r} \\
	\xi_3 - \xi_2 & \ldots & \frac{\xi_3^{r} - \xi_2^{r}}{r} \\
	\vdots & \ddots & \vdots \\
	\xi_{r+1} - \xi_{r} & \ldots & \frac{\xi_{r+1}^{r} - \xi_{r}^{r}}{r} \\	
\end{pmatrix} .$$ 
We may write the matrix $\widehat{V}^{\mathfrak{M}}(\mathcal{S}) $ now as the product 
\begin{equation*} \label{eq:detVDMnonorm}
\widehat{V}^{\mathfrak{M}}(\mathcal{S})  = \widehat{W}^{\mathfrak{M}}(\mathcal{S})  \mathbf{D}_{r},
\end{equation*}
with the two factors
$$ \widehat{W}^{\mathfrak{M}}(\mathcal{S}) =
\begin{pmatrix}
	\xi_2 - \xi_1 & \ldots & \xi_2^{r} -\xi_1^{r} \\
	\vdots & \ddots & \vdots \\
	\xi_{r+1} - \xi_{r} & \ldots & \xi_{r+1}^{r} - \xi_{r}^{r}
\end{pmatrix} \quad \text{and} \quad \mathbf{D}_{r} = \begin{pmatrix}
	1 & 0  & \ldots & 0 \\
	0 & \frac{1}{2}  & \ddots & 0 \\
	\vdots & \ddots  & \ddots & \vdots \\
	0 & 0  & \ldots & \frac{1}{r}
\end{pmatrix}. $$
We clearly have $ \det \mathbf{D}_{r} = \frac{1}{r!}$. By subtracting in $ \widehat{W} $ the $(r-1)$-th row from the $r$-th, the $(r-2)$-th from the $(r-1)$-th and so on, and then applying the Laplace Theorem with respect to the first column, we get that 
$$ \det \widehat{W}^{\mathfrak{M}}(\mathcal{S}) = \det V(\mathcal{X}) .$$ 
Hence, by Binet's Theorem and Eq. \eqref{eq:detVDMproduct}, we get
\begin{equation} \label{eq:detVdmsegments}
	\det \widehat{V}^{\mathfrak{M}}(\mathcal{S}) = \det \widehat{W}^{\mathfrak{M}}(\mathcal{S}) \cdot \det \mathbf{D}_{r} = \frac{1}{r!} \cdot \prod_{1 \leq i < j \leq r+1} (\xi_j - \xi_i).
\end{equation}
Since $ \frac{1}{r!}$ is a constant factor, we have thus proved the following result.

\begin{theorem}
The concatenated segments $\widehat{\mathcal{S}}^{\mathrm{Fek}}$ that solve the non-normalized Fekete Problem \ref{prob:concatenated} are uniquely determined as 
\[ s_i = [\xi_{i}^{(LL)}, \xi_{i+1}^{(LL)}], \quad i \in \{1, \ldots, r\},\]
where the points $\xi_{i}^{(LL)}$, $i \in \{1, \ldots, r+1\}$, denote the Legendre-Gauss-Lobatto nodes in $I$ that solve the nodal Fekete Problem \ref{prob:Fek} for point sets of size $r+1$.
\end{theorem}

\begin{figure}[H]
	\centering
	\includegraphics[width=8cm]{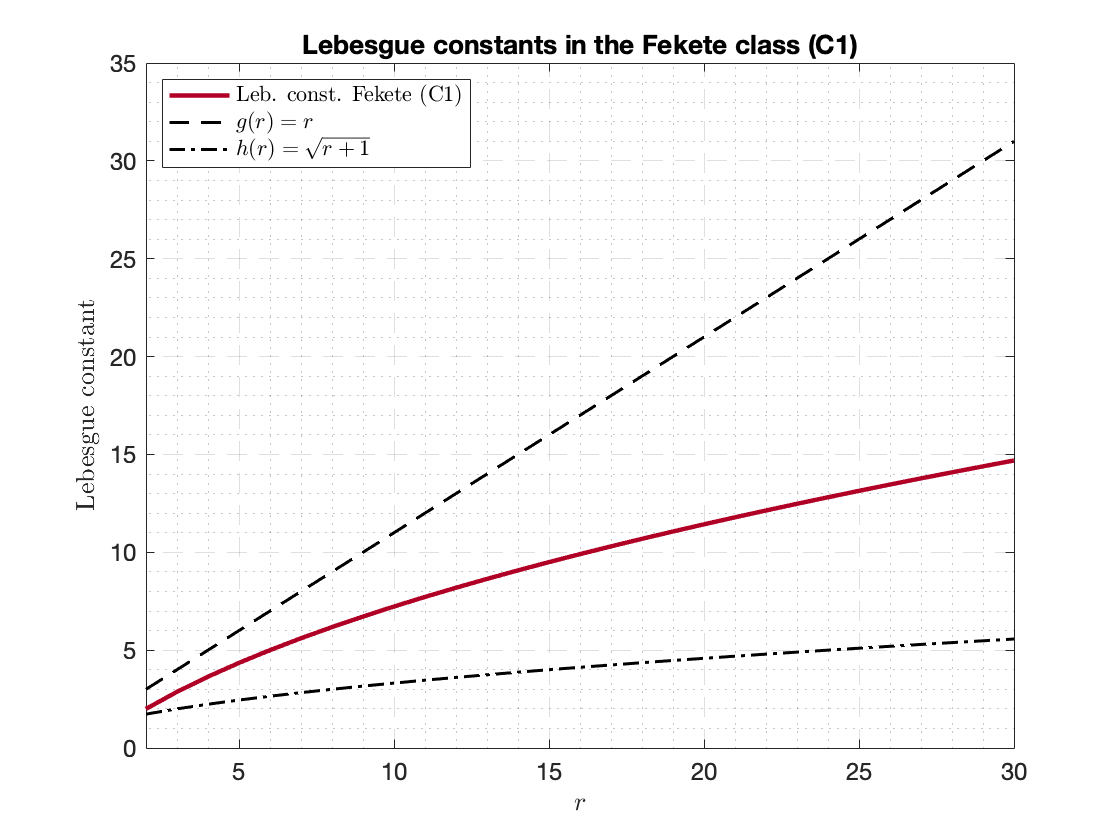}
	\includegraphics[width=8cm]{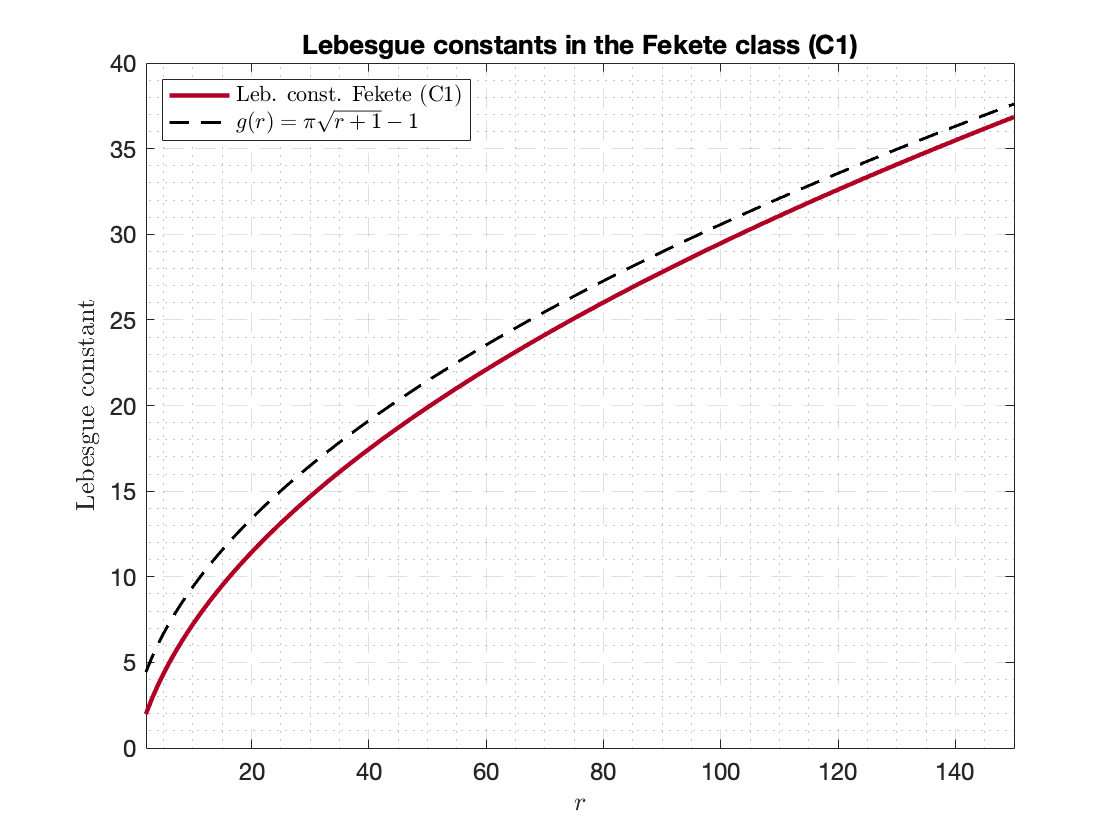}
	\caption{Lebesgue constants for concatenated (non-normalized) Fekete segments with degrees up to $ r = 150 $. The figures indicate a sublinear growth of the Lebesgue constant proportional to $\sqrt{r}$.}
	\label{fig:LebSeg}
\end{figure}

The growth of the Lebesgue constant \eqref{eq:segmentalLeb} for the Fekete segments $\widehat{\mathcal{S}}^{\mathrm{Fek}}$ is depicted in Figure \ref{fig:LebSeg}.
This growth behavior is consistent with the theoretical estimates 
$$ \Lambda_{r} (\widehat{\mathcal{S}}^{\mathrm{Fek}}) \leq r^3 \Lambda_{r+1} (\mathcal{X}^{\mathrm{Fek}})$$
of the Lebesgue constant for segments in the class (C1) derived in \cite[Proposition $5.1$]{BruniErb}. The comparison in Figure \ref{fig:LebSeg} indicates that the constant $\Lambda_{r} (\widehat{\mathcal{S}}^{\mathrm{Fek}})$ grows asymptotically as $\sqrt{r}$ rather than logarithmically. 
This is in contrast to the nodal setting in which the constant $\Lambda_{r+1} (\mathcal{X}^{\mathrm{Fek}})$ grows in fact logarithmically. In particular, this example shows that the Lebesgue constant for the end-points of segments might behave differently than the Lebesgue constant for the corresponding segments.

\subsection{The normalized Fekete problem for nodal-segmental interpolation}

For the general combined nodal-segmental interpolation we can formulate the normalized Fekete problem as follows.

\begin{problem} \label{prob:averaged}
	Find regular sets of supports $ \mathcal{A}^{\mathrm{Fek}}$ as introduced in Definition \ref{def:unisolvence} that maximise the determinant $|\det V(\mathcal{A})|$ of the normalized Vandermonde matrix in \eqref{eq:VandermondeAveraged}.
\end{problem}

As in the pure segmental and nodal framework, this problem is independent of the selected basis. In this general setting, there are at most $2r$ degrees of freedom that have to be handled in the optimization problem. For analytic studies, we restricted ourselves to scenarios in which the number of free parameters is reduced. We will treat as particularly interesting cases the special scenarios (C1) and (C2).  

\subsubsection{The case (C1): concatenated segments} \label{sect:FekCS}

\begin{figure}[H]
	\centering
	\includegraphics[width=8cm]{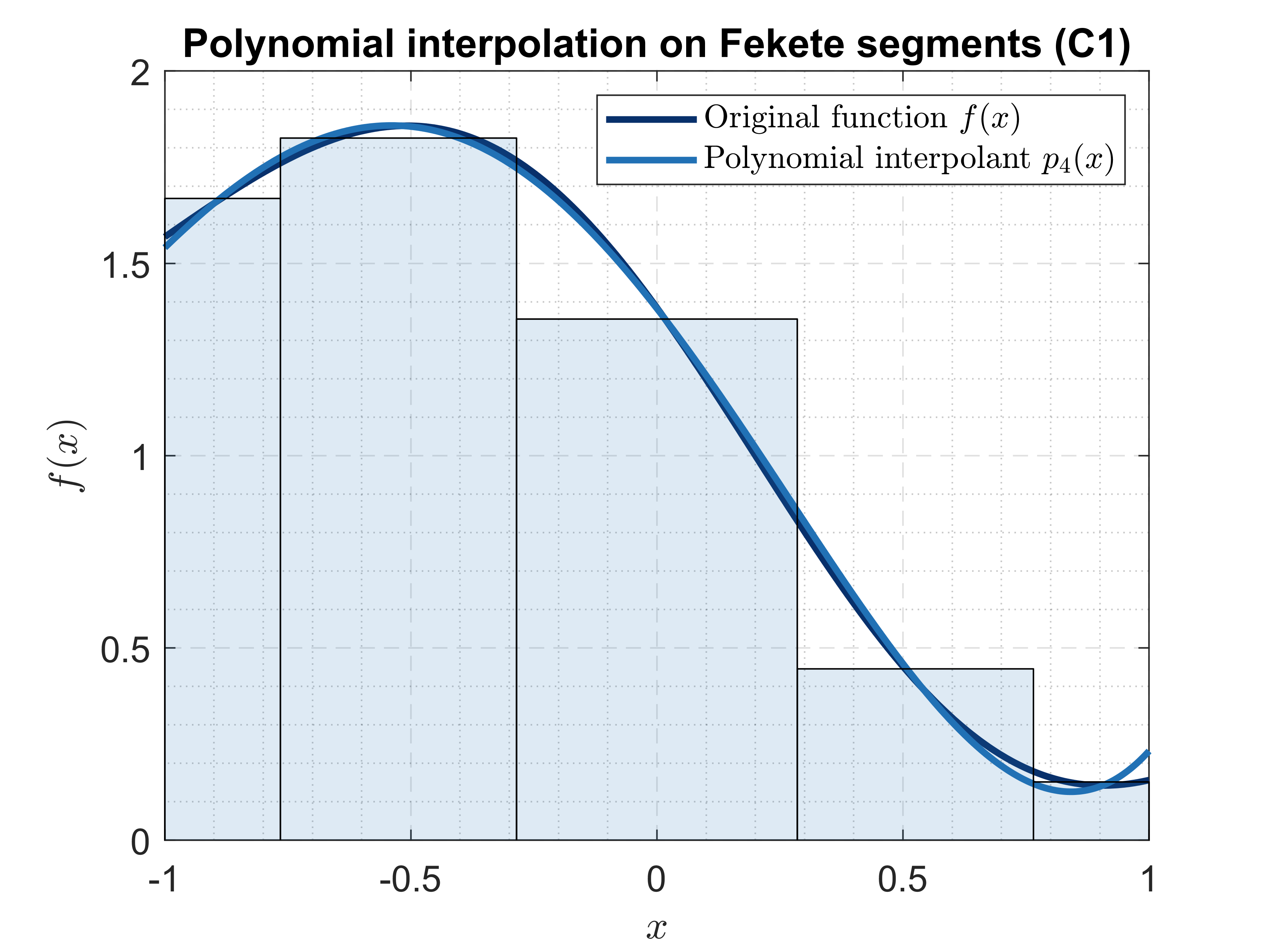}
	\includegraphics[width=8cm]{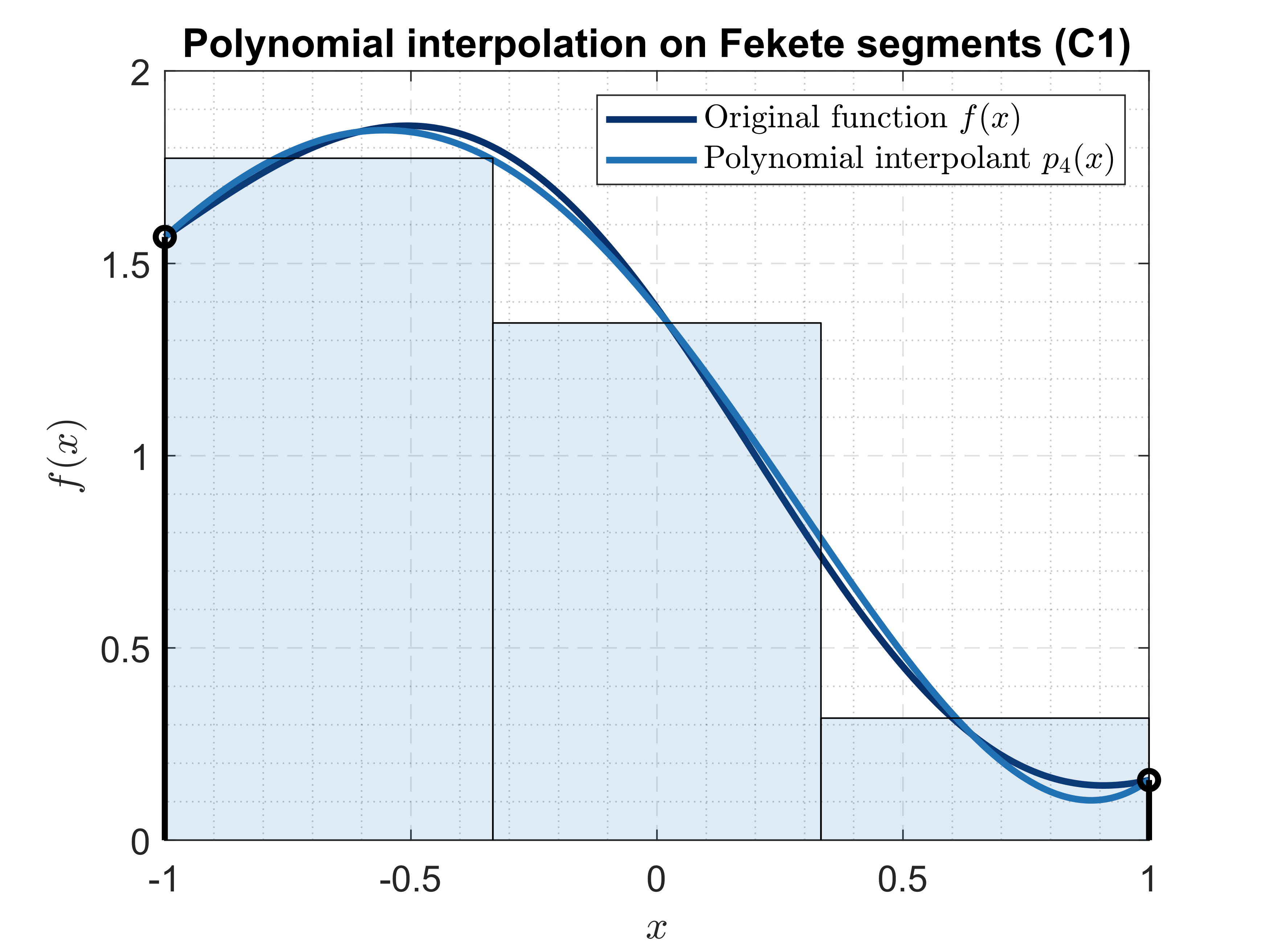}
	\caption{Segmental interpolation for Fekete segments in the class (C1) for the non-normalized Fekete (left) and the normalized Fekete problem (right). In the normalized case, the first and the last segment degenerate to single nodes. }
	\label{fig:FeketeC1}
\end{figure}

We consider again the monomial basis $\mathfrak{M} $ for $\P_{r-1}$ and recall that $ \mathbf{N} $ denotes the diagonal matrix with the segment lengths $|s_i|$ on the diagonal. Inserting concatenated segments $ s_i = [\xi_{i}, \xi_{i+1}] $ in the definition \eqref{eq:VandermondeAveraged} of the normalized Vandermonde matrix one retrieves that
\begin{equation*}
	V^{\mathfrak{M}}(\mathcal{S}) = \mathbf{N}^{-1} \widehat{V}^{\mathfrak{M}}(\mathcal{S}).
\end{equation*}
By Binet's theorem and \eqref{eq:detVdmsegments}, one thus finds that
\begin{align} \label{eq:detaveraged}
	\det V^{\mathfrak{M}} ({\mathcal{S}}) & = \det \mathbf{N}^{-1} \det \widehat{V}^{\mathfrak{M}}({\mathcal{S}}) = \frac{1}{\prod_{j=1}^{r} (\xi_{j+1}-\xi_{j})} \cdot \frac{1}{r!} \cdot \prod_{1 \leq i < j \leq r+1} (\xi_j - \xi_i) \\
	& = \frac{1}{r!} \cdot \prod_{1 \leq i + 1 < j \leq r+1} (\xi_j - \xi_i). \nonumber
\end{align}
The absence of the term $ i+1=j$ in this product makes the maximization of $|\det V^{\mathfrak{M}} ({\mathcal{S}})|$ very different from the maximization of the non-normalized determinant $|\det \widehat{V}^{\mathfrak{M}} ({\mathcal{S}})|$ in Problem \ref{prob:concatenated}. As a qualitative result, we get the following property of the concatenated Fekete segments $\mathcal{S}^{\mathrm{Fek}}$ of the normalized Fekete problem.

\begin{proposition} \label{prop:extremashrinking}
	If the Fekete segments $ \mathcal{S}^{\mathrm{Fek}} $ maximise the normalized Vandermonde determinant $|\det V^{\mathfrak{M}} ({\mathcal{S}})|$ given in \eqref{eq:detaveraged} for all concatenated segments $\mathcal{S}$ of size $r$, then $ \xi_1 = \xi_2 = -1 $ and $ \xi_{r} = \xi_{r+1} = 1 $, i.e., the first and the last element of the set $\mathcal{S}^{\mathrm{Fek}}$ are in fact single nodes and coincide with the endpoints of the interval $I$.  
\end{proposition}

{\setlength{\parindent}{0cm}
\textbf{Proof.}}
	Since $ \xi_j > \xi_i $ if $ j > i $, we have
	$$ \xi_{r} - \xi_j \leq \xi_{r+1} - \xi_j \leq 1 - \xi_j .$$
	Likewise,
	$$ \xi_j - \xi_2 \leq \xi_j - \xi_1 \leq \xi_j - (-1) = \xi_j + 1.$$
	Since the terms $(\xi_{r+1} - \xi_{r})$ and $(\xi_{2} - \xi_{1})$ are not contained in the product \eqref{eq:detaveraged}, the maximum of the product can only be attained if $\xi_{r+1} = \xi_{r} = 1$ and $\xi_{2} = \xi_{1} = -1$.   
\QED

The just proven proposition states that, in the solution of the normalized Fekete problem for concatenated segments, at least the two extremal segments collapse, i.e. $ \xi_0 = -1 $ and $ \xi_{r+1} = 1$ are the first and last element of $\mathcal{S}^{\mathrm{Fek}}$. A numerical simulation indicates that these are in fact the only two single nodes in $\mathcal{S}^{\mathrm{Fek}}$. Table \ref{tab:averagedFekete} additionally provides the explicit solutions in the concatenated case up to $r = 6$.
\begin{table}[!h]
	\begin{center}
		{\footnotesize
			\begin{tabular}{c|c|c|l}
				points & segments $r$ & polynomial degree & location of endpoints\\
				\hline
				2 & 1 & 0 & $\xi_1 = -1$, $ \xi_2 = 1$ \\
				3 & 2 & 1 & $\xi_1 = -1$, $ \xi_3 = 1$, independent of $ \xi_2 $ \\
				4 & 3 & 2 & $ \xi_1 = \xi_2 = -1 $, $ \xi_3 = \xi_4 = 1 $ \\
				5 & 4 & 3 & $ \xi_1 = \xi_2 = -1 $, $\xi_3 = 0$, $ \xi_4= \xi_5 = 1 $\\
				6 & 5 & 4 & $ \xi_1 = \xi_2 = -1 $, $ -\xi_3 = \xi_4 = \frac{1}{3} $, $ \xi_5 = \xi_6 = 1 $ \\
				7 & 6 & 5 & $ \xi_1 = \xi_2 = -1 $, $ -\xi_3 = \xi_5 = \frac{1+2\sqrt{2}}{7} $, $\xi_4 = 0 $, $ \xi_6 = \xi_7 = 1 $
		\end{tabular}}
		\caption{Endpoints of the segments $ s_i = [\xi_{i}, \xi_{i+1}]  $ in $ \mathcal{S}^{\mathrm{Fek}} $ that solve the Problem \ref{prob:averaged} under the restriction that the segments are concatenated, up to $ r = 6 $.} \label{tab:averagedFekete}
	\end{center}
\end{table}

Although we are not able to provide a closed formula for the segments that maximize the Vandermonde determinant \eqref{eq:detaveraged}, a solution of this Fekete problem can be calculated numerically. For this, we use an interior point method that computes the maximum $ \mathcal{S}^{\mathrm{Fek}}$ of the normalized Vandermonde determinant iteratively. This optimization problem is a non-linear constrained maximization problem. The Vandermonde determinant is not directly suited as a target function for the maximization procedure. Instead, we took the logarithm $\log \det V(\mathcal{S})$ as a target function together with the constraint that the endpoints $\xi_i$ of the segments have to be in the interior of $[-1,1]$. 

\begin{figure}[H]
	\centering
	\includegraphics[width=8cm]{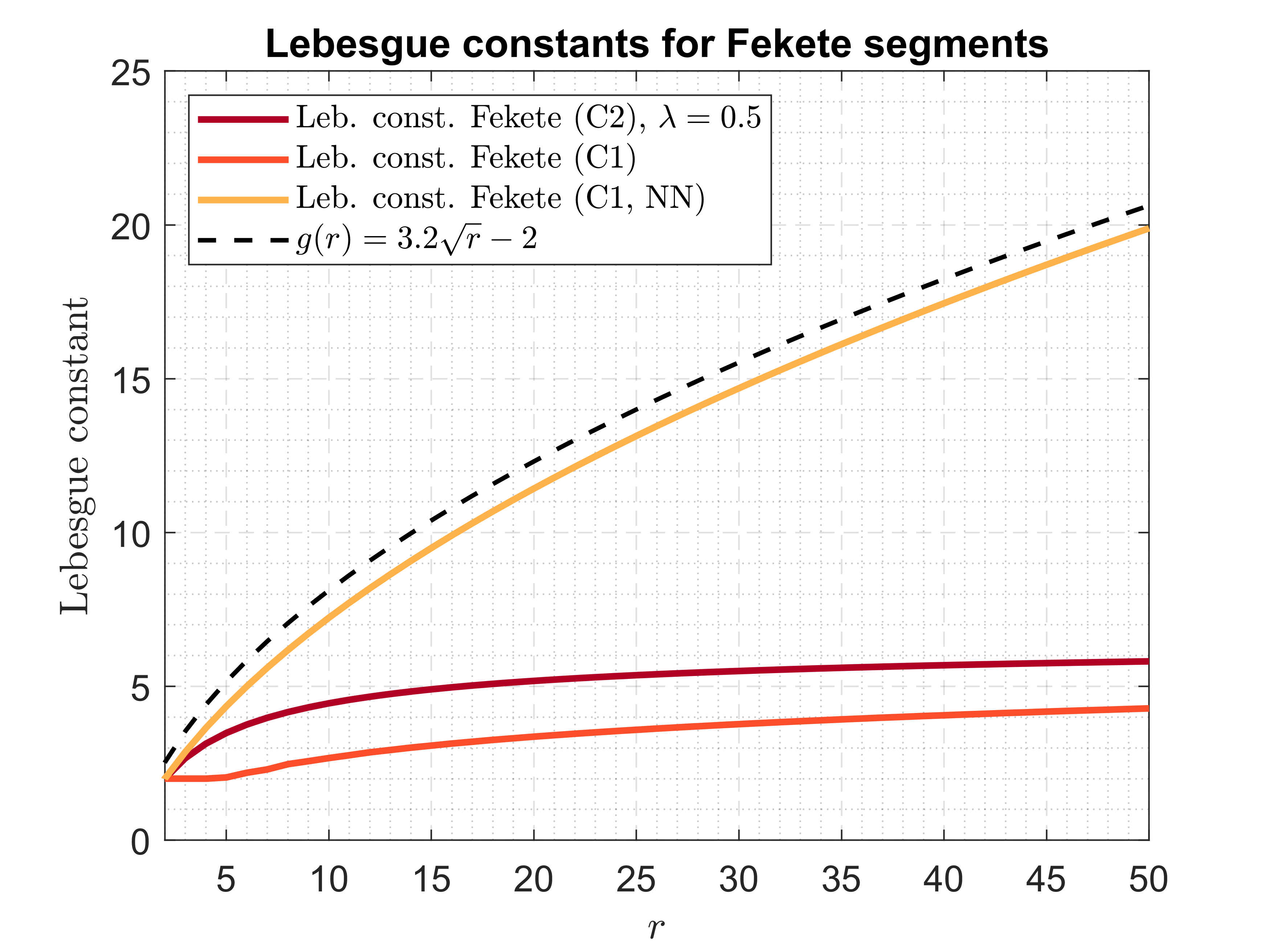}
	\includegraphics[width=8cm]{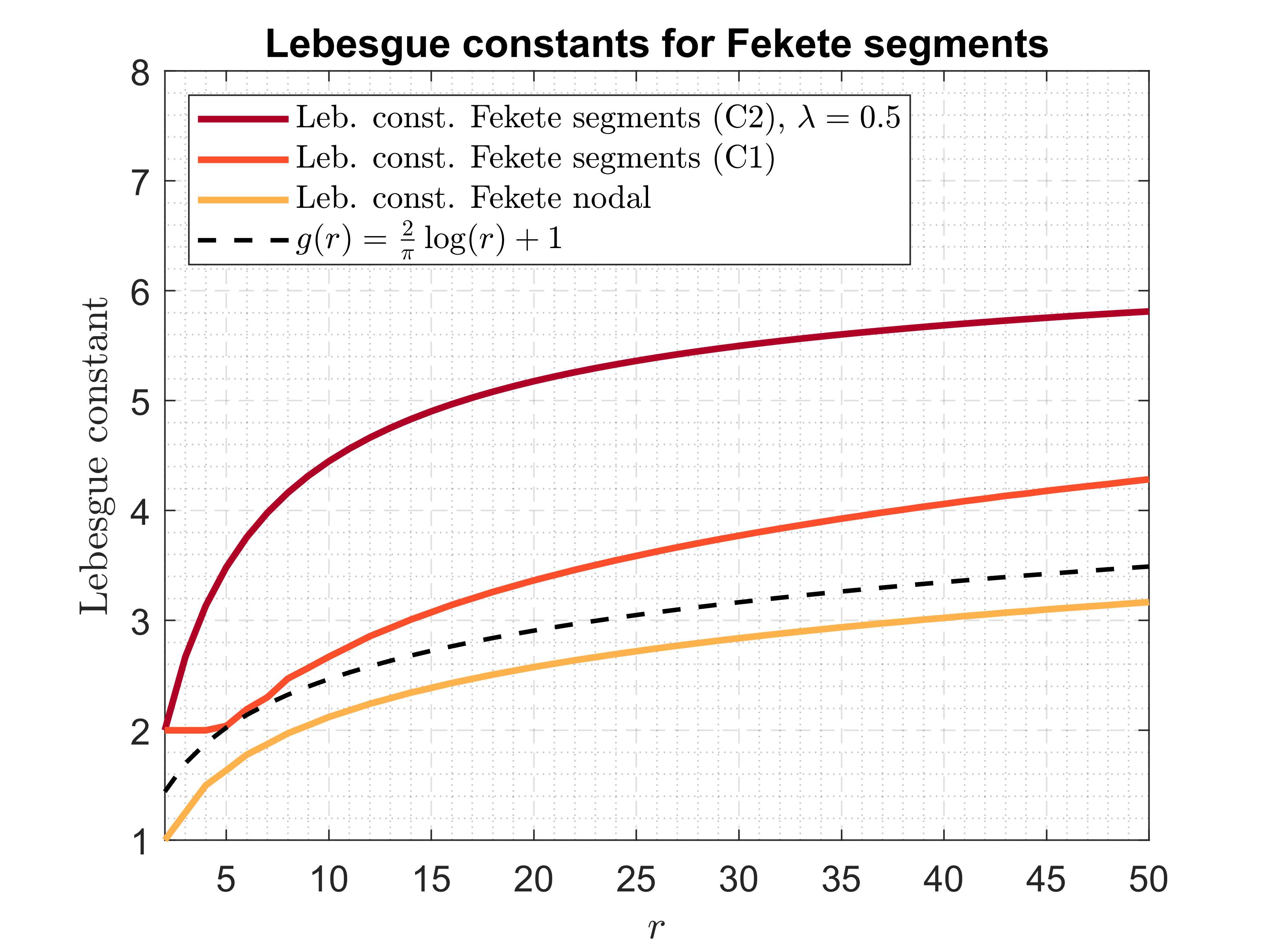}
	\caption{Comparison of the Lebesgue constants for Fekete segments in the classes (C1) and (C2). Left: comparison between normalized and non-normalized (NN) Fekete segments in the class (C1). Right: comparison between the Fekete segments in the class (C1) and the class (C2) for the values $\lambda = 0.5$ and $\lambda = 0$ (the nodal case).}
	\label{fig:LebesgueFeketeC1}
\end{figure}

The so-calculated numerical approximation of the Lebesgue constant is displayed in Figure \ref{fig:LebesgueFeketeC1}. It indicates that the Lebesgue constant for the normalized Fekete segments in the class (C1) displays a logarithmic growth, in contrast to the non-normalized case studied before.  

\subsubsection{The class (C2): segments with uniform arc-length} \label{sect:arc-length}

As a second relevant class of interval segments, we consider the class (C2) of segments with uniform arc-length. Similarly as in the scenario (C1), we are interested in resolving Problem \ref{prob:averaged} for the class (C2), taking the normalized Vandermonde matrix \eqref{eq:VandermondeAveraged} into account. This class is given as follows: for $r \in \N$, we consider a uniform arc-radius $0 < \rho < \pi/2$ (that may depend on $r$), and $r$ values $\rho \leq \tau_1 < \cdots < \tau_{r} \leq \pi-\rho$.
Based on these parameters, the endpoints of the interval segments $s_i$ in the class (C2) of segments with uniform arc-length are given as 
\begin{equation} \label{eq:C2} \alpha_i = \cos(\tau_i + \rho), \quad \beta_i = \cos(\tau_i - \rho)\qquad i \in \{1, \ldots, r\}. \end{equation}
The interval lengths of the segments $s_i = [\alpha_i, \beta_i]$ are not uniform since
\[ |s_i| = \cos(\tau_i - \rho) - \cos(\tau_i + \rho) = 2 \sin(\tau_i) \sin(\rho).\]
However, if the interval segments $s_i$ are mapped on the upper half part of the complex unit circle, then the arc-length of every mapped segment on the halfcircle is equal to $2 \rho$. This distance will be referred to as arc-length of the segments $s_i$. The values $\cos \tau_i$ will be called arc-midpoints. 

The restrictions $\rho \leq \tau_1$ and $\tau_{r} \leq \pi - \rho$ ensure that the arc-distance of the arc-midpoint $\cos \tau_i$  to one of the boundary points of $I$ is at least $\rho$ for all $i \in \{1, \ldots, r\}$. In the following, we will sometimes use the interval 
\[I_{\rho} = [-\cos \rho, \cos \rho ] \subset I,\]
and the latter restrictions are equivalent to the fact that all arc-midpoints lie in $I_{\rho}$.

For segments $\mathcal{S}$ in the class (C2), the basis $\mathfrak{U} = \{U_{j-1}\}_{j=1}^{r}$ of Chebyshev polynomials of the second kind defined by
\[U_{j-1}(x) = \frac{\sin( j \arccos x )}{ \sin( \arccos x)}, \quad j \in \N,\]
turns out to be advantageous for computational purposes, see \cite{BruniErb}. One reason is that for the Chebyshev basis $\mathfrak{U}$ and a set $\mathcal{S}$ of segments in the class (C2), the Vandermonde matrix $V^{\mathfrak{U}}(\mathcal{S})$ can be considerably simplified obtaining the entries
\begin{align} \label{eq:VandermMatrixChebyshev}
	V^{\mathfrak{U}}_{i,j} ({\mathcal{S}}) &=  \frac{1}{|a_i|}\int_{\alpha_i}^{\beta_i} U_j(x) \de x = \frac{\cos( j \arccos \alpha_i) - \cos( j \arccos \beta_i)}{j (\beta_i - \alpha_i)} = \frac{1}{j} \frac{\sin \left( j \tau_i \right)}{\sin \tau_i} \frac{\sin \left(j\rho \right)}{\sin \rho}.
\end{align}

In the class (C2), Fekete segments with uniform arc-length are uniquely characterized by the following result.

\begin{proposition} \label{prop-feketeC2}
	We consider the class (C2) of segments $\mathcal{S}$ of cardinality $|\mathcal{S}| = r $, $r \geq 2$, with uniform arc-radius $0< \rho < \pi/r$ such that all arc-midpoints $\cos \tau_i$ are in the subinterval $I_{\rho} = [-\cos \rho, \cos \rho]$. Then, the Fekete segments $s_i \in \safek $ in this class that maximize the determinant $| \det V^{\mathfrak{U}}({\mathcal{S}})|$ of the Vandermonde matrix, have the arc-midpoints
	\[\cos \tau_i^{\mathrm{Fek}} = \xi_{i}^{(LL)} \cos \rho, \quad i \in \{1, \ldots, r\},\]
	where $\xi_{i}^{(LL)}$, $i \in \{1, \ldots, r\}$ are the Legendre-Gauss-Lobatto nodes of order $r$ in $[-1,1]$. 
\end{proposition} 

{\setlength{\parindent}{0cm}
\textbf{Proof.}}
	In the light of the explicit expressions \eqref{eq:VandermMatrixChebyshev} and the definition of the Chebyshev polynomials $U_{j-1}$, we can factorize the Vandermonde matrix $V^{\mathfrak{U}} ({\mathcal{S}})$ as 
	\begin{equation} \label{eq:decompositionVMU} 
		V^{\mathfrak{U}} ({\mathcal{S}}) = 
		V^{\mathfrak{U}}(\mathcal{T}) \mathbf{D}_{\rho},  
	\end{equation}
	where $V^{\mathfrak{U}}(\mathcal{T})$ denotes the nodal Vandermonde matrix \eqref{eq:NodalVandermonde} with respect to the Chebyshev basis $\mathfrak{U}$ and the nodal evaluations at the arc-midpoints $\mathcal{T} = \{\cos \tau_1, \ldots, \cos \tau_{r}\}$, and $\mathbf{D}_{\rho}$ is a diagonal matrix given as
	\begin{equation} \label{eq:diagmatrixD}
		\mathbf{D}_{\rho} \doteq \begin{pmatrix} 1 &&& \\ & \frac{\sin \left(2\rho \right)}{2 \sin \rho} && \\ && \ddots & \\ &&& \frac{\sin \left(r\rho \right)}{r \sin \rho} \end{pmatrix}.
	\end{equation} 
	Since $\rho < \pi/r$, all diagonal entries of $\mathbf{D}_{\rho}$ are positive. 
	Therefore, as the arc-radius $\rho$ is fixed, the determinant in \eqref{eq:decompositionVMU} is maximized if and only if the determinant of the nodal Vandermonde matrix $V^{\mathfrak{U}}(\mathcal{T})$ is maximized. Combining Proposition \ref{lem:legendrezeros} with the affinity that maps $ [-1,1] $ to $ I_{\rho} $, we get that $|\det V^{\mathfrak{U}}(\mathcal{T})|$ is maximized for the Legendre-Gauss-Lobatto nodes of order $r$ scaled and shifted to the underlying interval which in our case is precisely the interval $I_{\rho}$. This provides the statement of the proposition.
\QED

\begin{figure}[h]
	\centering
	\includegraphics[width=8cm]{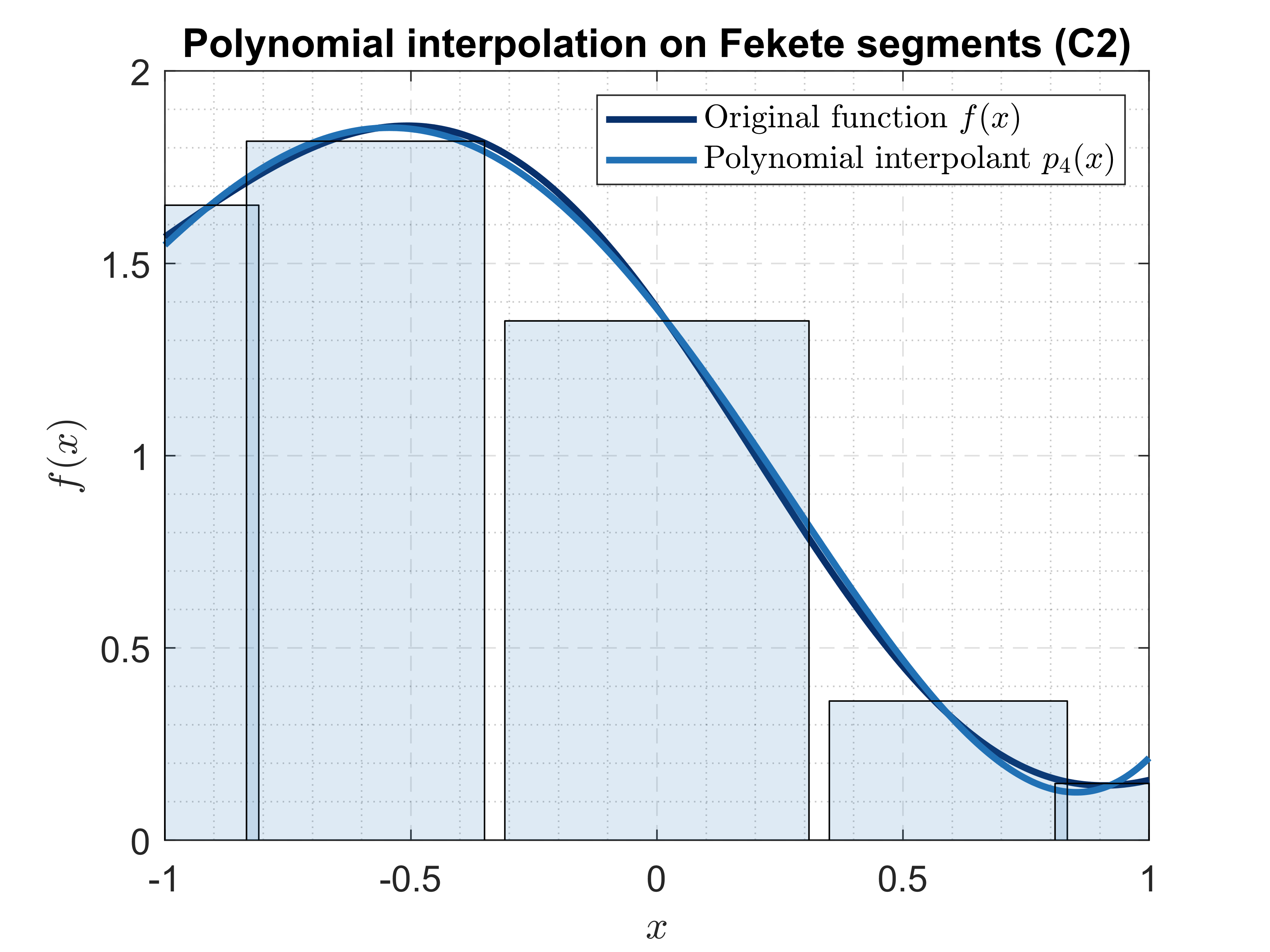}
	\includegraphics[width=8cm]{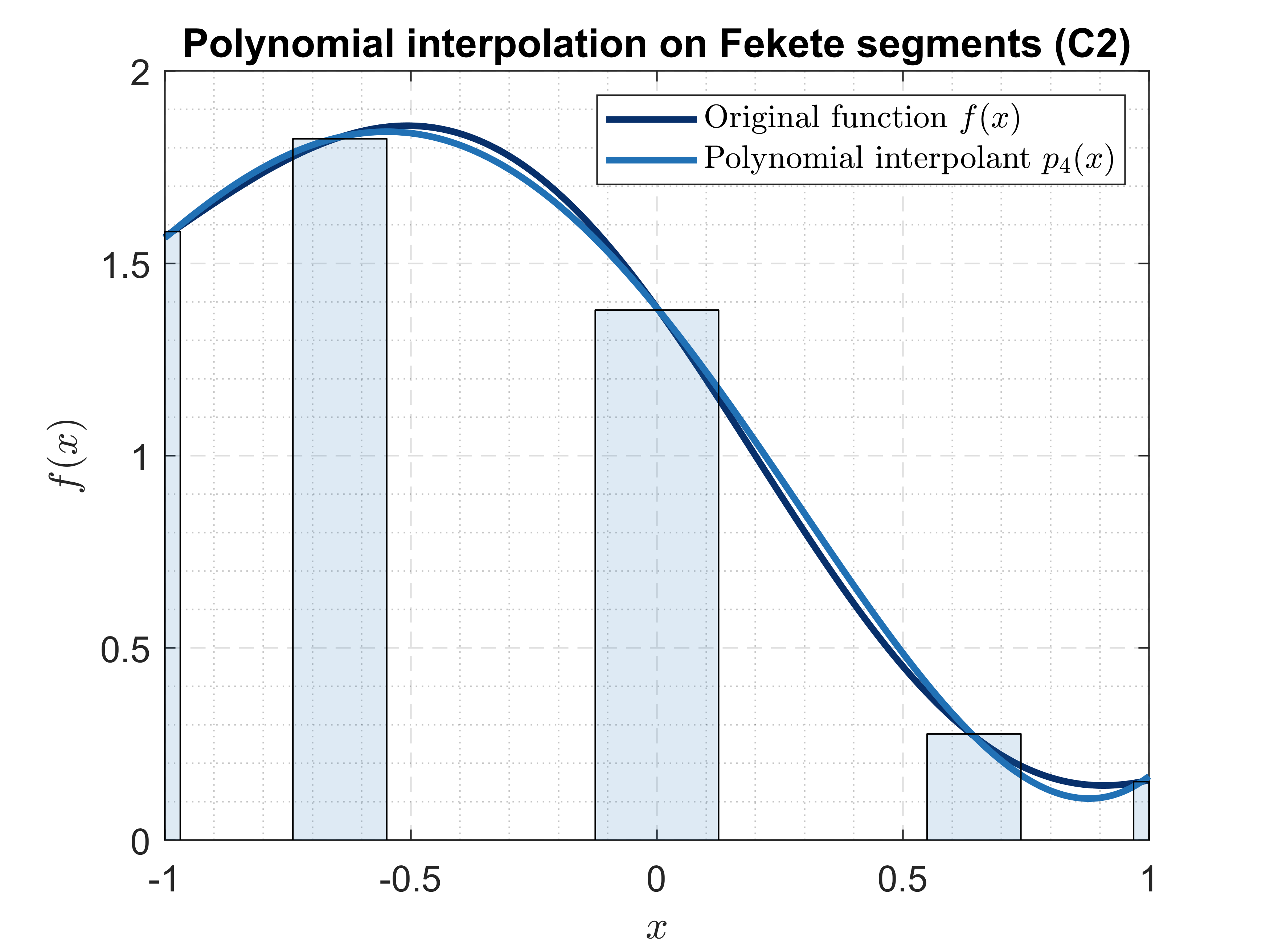}
	\caption{Interpolation for Fekete segments in the class (C2) with $\lambda = 0.5$ (left) and $\lambda = 0.2$ (right).}
	\label{fig:FeketeC2}
\end{figure}

For the solution of the polynomial interpolation problem, the Fekete segments $\safek$ in the class (C2) have very advantageous properties in terms of numerical conditioning, very similar to the classical Fekete nodes that correspond to the Legendre-Gauss-Lobatto nodes of order $r$ in $I$.
We can formalize this mathematically in terms of the norm
\[ \|\Pi_{r}(\safek)\|_{\mathrm{op}} = \sup_{\|f\| \leq 1} \sup_{x \in I} |\Pi_{r}(\safek)f(x)| \]
of the interpolation operator $\Pi_{r}(\safek)$ that maps a continuous function $f$ to its polynomial interpolant based on the averages of $f$ on the Fekete segments $s_i \in \safek$. If the Fekete segments $s_i$ are regular in the sense of Definition \ref{def:unisolvence}, i.e., the segments $s_i$ do at most intersect at their boundaries, we saw in Theorem \ref{thm:equaivalenceLebesgue} that this operator norm corresponds to the Lebesgue constant
\begin{equation*} \label{eq:genLeb}
	\Lambda_{r} (\safek) = \sup_{x \in I} \sum_{i = 1}^{r} |\ell_{s_i}(x)|. 
\end{equation*}

As an auxiliary tool to derive concrete estimates for the norm $\|\Pi_{r}(\safek)\|_{\mathrm{op}}$, we consider the integral operator $K_{\tau}: C(I) \to C(I)$ defined for $0 < \rho < \pi/2$ and $\cos t \in (-1,1)$ as
\begin{equation} \label{eq:integraloperator}
	K_{\rho} f(\cos t) \doteq \frac{1}{\cos (t - \rho) - \cos(t + \rho)} \int_{\cos(t + \rho)}^{\cos (t - \rho)}f(x) \de x. 
\end{equation}
The definition of the operator $K_{\rho}$ can be extended to the boundaries of $I$ in a straightforward way, by taking appropriate limits for the endpoints $\cos t = \pm 1$. Further, it is easy to check (see \cite{BruniErb}) that for the Chebyshev polynomials $U_j(x)$, $j \in \N_0$, one gets
\begin{align} \label{eq:spectraldecompositionK}  K_{\rho} U_j(\cos t) = \frac{1}{j+1}\frac{\sin((j+1)\rho)}{\sin(\rho)} U_j(\cos t), \notag
\end{align}
i.e., the polynomials $U_j$ are eigenfunctions of the integral operator $K_{\rho}$ with respect to the eigenvalues $\frac{1}{j+1}\frac{\sin((j+1)\rho)}{\sin(\rho)}$ and that $K_{\rho}$ maps the polynomial space $\mathbb{P}_{r-1}$ into $\mathbb{P}_{r-1}$ for every $r \in \N$. Further, the operator $K_{\rho}$ is invertible on $\mathbb{P}_{r-1}$ if $\rho < \pi/r$. We denote the restriction of $K_{\rho}$ to the polynomial space $\mathbb{P}_{r-1}$ by $K_{\rho,r} \doteq K_{\rho}\vert_{\mathbb{P}_{r-1}}$ and its respective inverse by $K_{\rho,r}^{-1}$. 

\vspace{1mm}

For the Fekete segments in the class (C2), we get the following quasi-optimal result on the asymptotic growth of the operator norm $\|\Pi_{r}(\safek)\|_{\mathrm{op}}$. 

\begin{proposition} \label{prop:LebesgueFeketeC2}
	Suppose $ r \geq 2 $. Let $\safek$ be the Fekete segments in the class (C2) with uniform arc-radius $\rho = \lambda \frac{\pi}{r}$, $0 <\lambda < 1$, and arc-midpoints $\mathcal{T}^{\mathrm{Fek}} = \{\cos \tau_1, \ldots, \cos \tau_{r}\} \subset I_{\rho}$. Then, the numerical conditioning of the interpolation problem $\| \Pi_{r}(\safek)\|_{\mathrm{op}}$ is bounded by
	\[ \frac1{2} \left( \frac{4}{\pi^2} \ln (r-1) - 1 \right) \leq 
	\| \Pi_{r}(\safek)\|_{\mathrm{op}} \leq  C_{\lambda} \, \ln (r)  , \]
	with a constant $C_{\lambda}$ that does not depend on the degree $r$.  
\end{proposition}

{\setlength{\parindent}{0cm}
\textbf{Proof.}}
	%Since with the given assumptions the interpolation operator $\Pi_{d+1}(\safek)$ is a projection from the space $C(I)$ onto the space $\mathbb{P}_d$, see Proposition \ref{prop:avgprojector}, a variant of the Kharshiladze-Lozinski theorem gives the lower inequality in Proposition \ref{prop:LebesgueFeketeC2} (see \cite[Chapter 3, Theorem 1]{CheneyLight2000}, the proof can be obtained from the derivations in \cite[Chapter 6, Section 5]{Cheney82}). 
	Proposition \ref{prop:avgprojector} ensures that the interpolation operator $\Pi_{r}(\safek)$ is a projection from the space $C(I)$ onto the space $\mathbb{P}_{r-1}$. Thus \eqref{eq:lowerboundLeb} gives the stated lower bound.
	
	We can thus proceed with the proof for the upper estimate. From the derivations given in \cite{BruniErb} we know that the interpolation operator $\Pi_{r}(\safek)$ for the Fekete segments $\safek$ in the class \eqref{eq:C2} can equivalently be written as
	\[ \Pi_{r}(\safek)f(x) = K_{\rho,r}^{-1} \Pi_{r}(\mathcal{T}^{\mathrm{Fek}}) K_{\rho}f(x)\]
	using the integral operator $K_{\rho}$ in \eqref{eq:integraloperator}, its inverse $K_{\rho,r}^{-1}$ for the polynomial space $\mathbb{P}_{r-1}$, and the nodal interpolation operator $\Pi_{r}(\mathcal{T}^{\mathrm{Fek}})$ on the arc-midpoints $\mathcal{T}^{\mathrm{Fek}}$ of the Fekete segments. The operator norm of $K_{\rho}$ on $C(I)$ is equal to $1$, while it is shown in \cite[Lemma 5.5]{BruniErb} that for $\rho = \frac{\lambda \pi}{r}$ with $0 < \lambda < 1$, the operator norm $\|K_{\rho,r}^{-1}\|_{\mathrm{op}}$ for the inverse is uniformly bounded and depends only on the parameter $\lambda$. We therefore obtain
	\[\| \Pi_{r}(\safek)\|_{\mathrm{op}} \leq \|K_{\rho,r}^{-1}\|_{\mathrm{op}} \| \Pi_{r}(\mathcal{T}^{\mathrm{Fek}})\|_{\mathrm{op}},\]
	and, thus, the upper estimate for the segmental operator norm boils down to the estimate of a nodal Lebesgue constant on the arc-midpoints $\mathcal{T}^{\mathrm{Fek}}$. 
	
	According to Proposition \ref{prop-feketeC2}, the arc-midpoints $\mathcal{T}^{\mathrm{Fek}}$ correspond to the Gauss-Legendre-Lobatto nodes on the interval $I_{\rho}$. For the respective Lebesgue constant, we know that (see \cite{Hesthaven,Sundermann})
	\[ \sup_{\|f\| \leq 1} \sup_{x \in I_{\rho}} |\Pi_{r}(\mathcal{T}^{\mathrm{Fek}})f (x) | \leq c \ln (r)\]
    with a finite but unknown constant $c > 0$.	
	%The constant $c > 0$ in this inequality is finite but unknown. Numerical computations performed in \cite{Hesthaven} however indicate that the resulting Lebesgue constant is in general smaller than the Lebesgue constant for the Chebyshev-Lobatto nodes. 
	Finally, we have to extend the inequality above from the subinterval $I_{\rho}$ to the full interval $I$. This can be done by using the fact that the Chebyshev polynomials $T_{r-1}$ of the first kind scaled to the interval $I_{\rho}$ are the polynomials of degree $r-1$ that grow most rapidly outside $I_{\rho}$ \cite[Section 2.7.1]{Rivlin}. Using this extremal property of the Chebyshev polynomials, we get for every polynomial $p(x)$ of degree $r-1$ the inequality
	\[ \sup_{x \in I} |p(x)| \leq \sup_{x \in I_{\rho}} |p(x)| \sup_{x \in I} |T_{r-1}\left(\frac{x}{\cos \rho}\right)| \leq \sup_{x \in I_{\rho}} |p(x)| \;  \left|T_{r-1}\left(\frac{1}{\cos \rho}\right)\right|.\]
	Using basic estimates for the Chebyshev polynomials outside the interval $I = [-1,1]$, we get
	\begin{align*}
		\left|T_{r-1}\left(\textstyle \frac{1}{\cos \rho}\right)\right| &= \cosh \left( (r-1) \;\mathrm{arccosh} \left(\textstyle \frac{1}{\cos \rho}\right) \right)  
		= \cosh \left( (r-1) \;\ln \left(\textstyle \frac{1 + \sin \rho}{\cos \rho}\right) \right) \\
		&\leq \cosh \left( (r-1) \;\ln \left( 1 + \sin \rho \right) \right) \leq \cosh \left( \;\ln \left( 1 + \frac{\lambda \pi}{r} \right)^{r-1} \right) \\
		& \leq \cosh \left(\ln \exp(\lambda \pi) \right) = \cosh (\lambda \pi). 
	\end{align*} 
	We therefore get the bound $ \|\Pi_{r}(\mathcal{T}^{\mathrm{Fek}})\|_{\mathrm{op}} \leq c \ln (r) \cosh (\lambda \pi)$ for the nodal Lebesgue constant on the arc-midpoints $\mathcal{T}^{\mathrm{Fek}}$, and thus
	\[ \| \Pi_{r}(\safek)\|_{\mathrm{op}} \leq c \cosh (\lambda \pi) \|K_{\rho,r}^{-1}\|_{\mathrm{op}}  \ln (r) \]
	for the segmental interpolation operator on the Fekete segments. Herein, we get a uniform estimate for the operator norm $\|K_{\rho,r}^{-1}\|_{\mathrm{op}}$ that depends only on $\lambda$ but not on $r$. 
\QED

For the class (C2), the Fekete segments can additionally be characterized in the following alternative way due to Fej\'er \cite{Fejer1932}. 

\begin{proposition}
	Let $\mathcal{S} = \{s_1, \ldots, s_{r}\}$ denote a set of $r$ segments in the class (C2) with uniform arc-radius $0 < \rho < \frac{\pi}{r}$, and arc-midpoints $\mathcal{T} = \{\cos \tau_1, \ldots, \cos \tau_{r}\} \subset I_{\rho}$. 
	Then,
	$$
	\min_{\mathcal{S}} \max_{x \in I_{\rho}}\left\{\left(K_{\rho} \ell_{s_1}(x)\right)^2 + \ldots +\left(K_{\rho} \ell_{s_{r}}(x)\right)^2\right\} = 1 .
	$$
	The only set of segments $\mathcal{S}$ for which the minimum is attained, i.e., for which
	\begin{equation} \label{eq:Fejersegments}
	\max_{x \in I_{\rho}} \left\{ \left(K_{\rho} \ell_{s_1}(x)\right)^2 + \ldots +\left(K_{\rho} \ell_{s_{r}}(x)\right)^2 \right\} = 1
	\end{equation}
	holds true, is the set $\safek$ of Fekete segments in the class (C2). In the limit $\rho \to 0$, this corresponds to the classical characterization of the Fekete nodes in the interval $I$. 
\end{proposition}

{\setlength{\parindent}{0cm}
\textbf{Proof.}}
	If the segments $s_i$ are in the class (C2), the Lagrange basis polynomials $\ell_{s_i}(x)$, $i \in \{1,\ldots,r\}$, satisfy the identity
	\[K_{\rho} \ell_{s_i}(x) = \ell_{\cos \tau_i}(x),\]
	i.e., the operator $K_{\rho}$ applied to the segmental Lagrange polynomials $\ell_{s_i}$ yields the nodal Lagrange polynomials $\ell_{\cos \tau_i}$ for the set $\mathcal{T}$ of arc-midpoints. This holds particularly true if the segments $s_i$ are the Fekete segments in $\safek$. On the other hand, it is well-known (see \cite{Fejer1932}) that the nodal Fekete nodes $\mathcal{T}^{\mathrm{Fek}}$ on the set $I_{\rho}$ are the only nodes such that  
	$$
	\min_{\mathcal{T} \subset I_{\rho}} \max_{x \in I_{\rho}}\left\{\ell_{\cos \tau_1}(x)^2 + \ldots + \ell_{\cos \tau_{r}}(x)^2\right\} = 1
	$$
	is satisfied.
\QED

\subsubsection{The class (C2) and a free parameter $\rho$}

Finally, in addition to the arc-midpoints in the class (C2), we allow also the arc-radius $\rho$ to be a free parameter in the Fekete maximization problem, and look for segments that maximize the normalized Vandermonde determinant over this larger class. It turns out that in this case the optimum is attained in the limit $\rho \to 0$ with
the Legendre-Gauss-Lobatto nodes $\xi_{i}^{(LL)}$ of order $r$ providing the solution of the Fekete problem. 

\begin{proposition}
	We consider the Legendre-Gauss-Lobatto nodes $\xi_{i}^{(LL)}$, $i \in \{1, \ldots, r\}$ (corresponding to the Fekete nodes in $I$) as limiting case of the Fekete segments $\safek$ in the class (C2) when the arc-radius $\rho \to 0$ tends to zero. 
	
	In the class \eqref{eq:C2} with $0 \leq \rho < \pi/r$, the Vandermonde determinant 
	$\det V(\safek)$ is a decreasing function of the arc-radius $\rho$. The maximal Vandermonde determinant among all amissible values $0 \leq \rho < \pi/r$ is thus obtained in the limiting case $\rho = 0$, i.e. when the Fekete segments degenerate to the Fekete nodes in $I$.  
\end{proposition}

\begin{remark}
	A similar monotonic behaviour can be observed numerically for the Lebesgue constant 
	$\Lambda_r (\safek)$. In this case, the values $\Lambda_r (\safek)$ increase when $\rho$ gets larger with a minimal value obtained for the nodal Lebesgue constant at the Fekete nodes in $I$.   
\end{remark}

{\setlength{\parindent}{0cm}
\textbf{Proof.}}
	As in the proof of Proposition \ref{prop-feketeC2}, we can use the Chebyshev basis $\{U_{j-1}\}_{j = 1}^{r}$ and factorize the Vandermonde matrix $V^{\mathfrak{U}} (\safek) = V^{\mathfrak{U}}(\mathcal{T}^{\mathrm{Fek}}) \mathbf{D}_{\rho} $, where $V^{\mathfrak{U}}(\mathcal{T}^{\mathrm{Fek}})$ denotes the Vandermonde matrix with respect to the Chebyshev basis $\mathfrak{U}$ and the nodal evaluations at the arc-midpoints $\mathcal{T}^{\mathrm{Fek}} = \{\cos \tau_1, \ldots, \cos \tau_{r}\}$ corresponding to the Fekete nodes in the subinterval $I_{\rho}$. This leads to the factorization
	$$ \det V^{\mathfrak{U}} (\safek) = \det V^{\mathfrak{U}}(\mathcal{T}^{\mathrm{Fek}}) \det \mathbf{D}_{\rho}. $$
	From the definition of the diagonal matrix $\mathbf{D}_{\rho}$ in \eqref{eq:diagmatrixD} we see that the determinant $\det \mathbf{D}_{\rho} = \prod_{i = 1}^{r} \frac{\sin \left(i \rho \right)}{i \sin \rho}$ is a decreasing function for increasing values $\rho < \pi/r$. Also, if $0\leq \rho < \rho' < \pi/2$, the interval $I_{\rho'}$ is contained in $I_{\rho}$ and we know that the Fekete nodes in $I_{\rho}$ maximize the determinant $\det V^{\mathfrak{U}}(\mathcal{T})$ among all node sets $\mathcal{T}$ in $I_{\rho}$. Therefore, also the second factor $\det V^{\mathfrak{U}}(\mathcal{T}^{\mathrm{Fek}})$ is a decreasing function of the arc-radius $\rho$. This provides the stated result. 
\QED

\section{Conclusions}

Motivated by the classical nodal Fekete problem, we studied corresponding problems in a generalized framework where input data is given by function averages over segments or by a combination of nodal and segmental data. In order to obtain families of supports that are valuable for univariate polynomial interpolation we were able to solve particular segmental Fekete problems explicitly. 
In the segmental setting we found out that a normalisation of the data by the measure of the supports is essential in retrieving the logarithmic behaviour of the Lebesgue constant. We constructed further a set of supports that contains an interesting warning: also if the end-points of the segments are distributed in such a way that the nodal Lebesgue constant increases slowly, the Lebesgue constant for the segmental problem may increase considerably faster. This warning should in particular be taken into account in higher-dimensional settings where supports are often constructed by connecting well-behaved sets of vertices. 

\section*{Acknowledgements}

This research has been accomplished within the research networks RITA and UMI-TAA, and was partially funded by GNCS-IN$\delta$AM. The first author is funded by IN$\delta$AM and supported by Università di Padova. The second author was funded by the European Union - NextGenerationEU under the National Recovery and Resilience Plan (NRRP), Mission 4 Component 2 Investment 1.1 - Call PRIN 2022 No. 104 of February 2, 2022 of Italian Ministry of University and Research; Project 2022FHCNY3 (subject area: PE - Physical Sciences and Engineering) "Computational mEthods for Medical Imaging (CEMI)"


\begin{thebibliography}{00}
	
	%\bibitem{Achieser92}
	%N.~I. Achieser,
	%\newblock \emph{Theory of Approximation}, 
	%\newblock Dover Publications, New York, 1992.
	
	%\bibitem{platte}
	%B.~Adcock and R.~B. Platte,
	%\newblock {\em A mapped polynomial method for high-accuracy approximations %on arbitrary grids}, 
	%\newblock {SIAM J. Numer. Anal.}, \textbf{54} (2016), pp. 2256--2281.
	
	\bibitem{BruniEdges}
	A.~Alonso~Rodr{\'i}guez, L.~Bruni~Bruno and F.~Rapetti,
	\newblock {\em Towards nonuniform distributions of unisolvent weights for {W}hitney finite element spaces on simplices: the edge element case},
	\newblock {Calcolo}, {\bf 59}(4):37 (2022).
	
	\bibitem{BruniRunge}
	A. Alonso Rodr\'iguez, L.~Bruni~Bruno and F. Rapetti,
	\newblock {\em Whitney edge elements and the Runge phenomenon},
	\newblock {J. Comput. Appl. Math.}, {\bf 427}:115117 (2023).
	
%	\bibitem{BruniFaces}
%	A. Alonso Rodr\'iguez, L.~Bruni~Bruno and F. Rapetti,
%	\newblock {\em Flexible weights for high order face based finite element interpolation},
%	\newblock in Spectral and High Order Methods for Partial Differential Equations ICOSAHOM 2020+1, J. M. Melenk, I. Perugia, J. Schöberl and C. Schwab, eds., Lect. Notes Comput. Sci. Eng. {\bf 137}, Spinger, Cham., 2023, pp. 117--128.
	
	\bibitem{AnaFra}
	A.~Alonso~Rodr{\'i}guez and F.~Rapetti,
	\newblock {\em On a generalization of the {L}ebesgue's constant},
	\newblock {J. Comput. Phys.}, {\bf 428}:109964 (2021).
	
    \bibitem{Bojanov}
    B. Bojanov,
    \newblock {Interpolation and integration based on averaged values},
    \newblock {in Approximation and probability}, Polish Acad. Sci. Inst. Math., Warsaw, {\bf 72} (2006), pp. 25--47.
%	
%	\bibitem{Bos1990}
%	L. Bos,
%	\newblock {Some remarks on the {F}ej\'{e}r problem for {L}agrange interpolation in several variables},
%	\newblock {J. Approx. Theory}, {\bf 60}(2) (1990), pp. 133--140.
%	
%	\bibitem{BosAlgebraic}
%	L. Bos,
%	\newblock {On certain configurations of points in {${\bf R}^n$} which are unisolvent for polynomial interpolation},
%	\newblock {J. Approx. Theory}, {\bf 64}(3) (1991), pp. 271--280.

    \bibitem{BosSurvey}
    L. P. Bos,
    \newblock {\em On Fekete points for a real simplex},
    \newblock {Indag. Math. (N.S.)}, {\bf 34}(2) (2023), pp. 274--293

    \bibitem{Bos2010}
	L. P. Bos, S. De Marchi, A. Sommariva and M. Vianello,
	\newblock {\em Computing multivariate Fekete and Leja points
     by numerical linear algebra},
	\newblock {SIAM J. Num. Anal.}, {\bf 48}(5) (2010), pp. 1984--1999.

	\bibitem{Bos1D}
	L. P. Bos and N. Levenberg,
	\newblock {\em On the calculation of approximate {F}ekete points: the univariate case},
	\newblock {Electron. Trans. Numer. Anal.}, {\bf 30} (2008), pp. 377--397.

 	\bibitem{Bos2001}
	L. P. Bos, M. Taylor and B. Wingate,
	\newblock {\em Tensor product Gauss-Lobatto points are Fekete points for the cube},
	\newblock {Math. Comput.}, {\bf 70} (2001), pp. 1543--1547.

	\bibitem{Bossavitbook}
	A. Bossavit,
	\newblock{Computational electromagnetism}. 
	\newblock{Academic Press, Inc., San Diego, CA}, 1998
	
	%\bibitem{AFW06}
	%D.~N.~Arnold, R.~S.~Falk and R.~Winther,
	%\newblock \emph{Finite Element exterior calculus, homological techniques, %and applications},
	%\newblock Acta Numer. \textbf{15} (2006), pp. 1--155.
	
	%\bibitem{Atkinson1998}
	%K. Atkinson and W. Han,
	%\newblock \emph{Theoretical Numerical Analysis: A Functional Analysis %Framework},
	%\newblock Springer, New York, 2009.
	
	%\bibitem{BossavitBook}
	%A.~Bossavit,
	%\newblock {\em Computational electromagnetism},
	%\newblock {Electromagnetism. Academic Press, Inc., San Diego, CA, 1998.}
	
	\bibitem{BruniThesis}
	L.~Bruni~Bruno,
	\newblock {\em Weights as degrees of freedoom for high order Whitney finite elements},
	\newblock Ph.D. thesis, University of Trento, 2022.
	
	%\bibitem{ComputingWeights}
	%L.~Bruni~Bruno, A.~Alonso~Rodr\'{\i}guez and F.~Rapetti,
	%\newblock {\em Computing weights for high order {W}hitney edge elements},
	%\newblock {Dolomites Res. Notes Approx.}, {\bf 15} (2022), pp. 1--12.
	
	\bibitem{BruniErb}
	L.~Bruni~Bruno and W. Erb,
	\newblock {\em Polynomial interpolation of function averages on interval segments},
	\newblock submitted (2023). ArXiv: \url{https://arxiv.org/abs/2309.00328}.
	
	%	\bibitem{BBSSC}
	%	L.~Bruni~Bruno, M. Semplice and S. Serra-Capizzano,
	%	\newblock {\em The numerical linear algebra of weights: from the spectral analysis to conditioning and preconditioning in the Laplacian case},
	%	\newblock submitted (2023). ArXiv: \url{https://arxiv.org/abs/2311.01467}.
	
	%\bibitem{Brutman1}
	%L.~Brutman,
	%\newblock {\em Lebesgue functions for polynomial interpolation -- a survey},
	%\newblock {Ann. Numer. Math.}, {\bf 4} (1997), pp. 111--127.
	
	%\bibitem{Cappellazzo2023}
	%G.~Cappellazzo, W.~Erb, F.~Marchetti and D.~Poggiali,
	%\newblock {\em On Kosloff Tal-Ezer least-squares quadrature formulas},
	%\newblock {BIT Numer. Math.}, {\bf 63}:15 (2023).
	
%	\bibitem{ChenBabuska}
%	Q. Chen and I. Babuška,
%	\newblock {Approximate optimal points for polynomial interpolation of real functions in an interval and in a triangle}, Comput. Methods Appl. Mech. Engrg. {\bf 128}:3-4 (1995), pp. 405--417.
	
	\bibitem{Cheney82}
	E. W. Cheney,
	\newblock \emph{Introduction to Approximation Theory},
	\newblock AMS Chelsea Publishing, Providence, RI, 1982.
	
	\bibitem{CheneyLight2000}
	E. W. Cheney and W. Light,
	\newblock \emph{A Course in Approximation Theory},
	\newblock AMS, Providence, RI, 2000.
	
	\bibitem{Davis75}
	P. J. Davis,
	\newblock \emph{Interpolation and Approximation},
	\newblock Dover Publications, New York, 1975.
	
	%\bibitem{DeMarchi2020}
	%S.~De~Marchi, F.~Marchetti, E.~Perracchione and D.~Poggiali,
	%\newblock {\em Polynomial interpolation via mapped bases without %resampling},
	%\newblock {J. Comput. Appl. Math.}, {\bf 364}:112347 (2020).

 	\bibitem{Fejer1932}
	L. Fej{\'e}r,
	\newblock{\em {B}estimmung derjenigen {A}bszissen eines {I}ntervalles f\"ur
welche die {Q}uadratsumme der {G}rundfunktionen der {L}agrangeschen
{I}nterpolation im {I}ntervalle $[-1,1]$ ein m\"oglichst kleines Maximum besitzt},
	\newblock {Ann. Scula Norm. Sup. Pisa Sci. Fis. Mat. Ser. II.}, {\bf 1} (1932), pp. 263--276.
	
	\bibitem{Fejer1932b}
	L. Fej{\'e}r,
	\newblock{\em Lagrangesche Interpolation und die Zugeh{\"o}rigen {K}onjugierten {P}unkte},
	\newblock {Math. Ann.}, {\bf 106} (1932), pp. 1--55.
	
	\bibitem{Fekete}
	M. Fekete,
	\newblock{\em \"{U}ber die {V}erteilung der {W}urzeln bei gewissen
		algebraischen {G}leichungen mit ganzzahligen {K}oeffizienten},
	\newblock {Math. Z.}, {\bf 17}(1) (1923), pp. 228--249
	
	\bibitem{Gerritsma}
	M. Gerritsma,
	\newblock {\em Edge functions for spectral element methods},
	\newblock in Spectral and high order methods for partial differential equations, Selected papers from the ICOSAHOM '09 conference, J. S. Hesthaven and E. M. R{\o}nquist, eds., Lect. Notes Comput. Sci. Eng. {\bf 76}, Spinger, Heidelberg, 2011, pp. 199--207.
	
	\bibitem{Harrison}
	J. Harrison,
	\newblock {\em Continuity of the integral as a function of the domain},
	\newblock {J. Geom. Anal.}, {\bf 8}:5 (1998), pp. 769--795.
	
	\bibitem{Hesthaven}
	J.~S.~Hesthaven,
	\newblock {\em From electrostatics to almost optimal nodal sets for polynomial interpolation in a simplex},
	\newblock {SIAM J. Numer. Anal.}, {\bf 35} (1998), pp. 655--676.
	
	\bibitem{Hiemstra2014}
	R. R. Hiemstra, D. Toshniwal, R. H. M. Huijsmans and M. I. Gerritsma,
	\newblock {\em High order geometric methods with exact conservation properties},
	\newblock {J. Comput. Phys.} {\bf 257} (2014), part B, 1444--1471.
	
	%\bibitem{Hiptmair2007}
	%R.~Hiptmair and J.~Xu,
	%\newblock {\em Nodal auxiliary space %preconditioning in H (curl) and H (div) spaces},
	%\newblock {SIAM J. Numer. Anal.}, {\bf 45}(6) (2007), pp. 2483-2509.
	
%	\bibitem{Pazner2023}
%	W.~Pazner, T.~Kolev and C.~R.~Dohrmann,
%	\newblock {\em Low-order preconditioning for the high-order finite element de Rham complex},
%	\newblock {SIAM J. Sci. Comput.}, {\bf 45}(2) (2023), pp. A675--A702.
%	
	\bibitem{Ibrahimoglu}
	B. A. Ibrahimoglu,
	\newblock {Lebesgue functions and {L}ebesgue constants in polynomial
		interpolation}
	\newblock{J. Inequal. Appl.}, (2016), Paper No. 93, pp. 15.
%	
%	\bibitem{m2an}
%	F.~Rapetti,
%	\newblock {\em High order edge elements on simplicial meshes},
%	\newblock {ESAIM Math. Model. Numer. Anal.}, {\bf 41}(6) (2007), pp. 1001--1020.
	
	%\bibitem{RB09}
	%F.~Rapetti and A.~Bossavit,
	%\newblock {\em Whitney forms of higher degree},
	%\newblock {SIAM J. Numer. Anal.}, {\bf 47}(3) (2009), pp. 2369--2386.
	
	%   \bibitem{Runge}
	%   C.~Runge,
	%   \newblock \emph{\"Uber empirische Funktionen und die Interpolation %zwischen \"aquidistanten Ordinaten},
	%   \newblock Zeit. Math. Phys. \textbf{46} (1901), pp. 224--243.
	
	%\bibitem{Schoenage}
	%A.~Sch\"{o}nhage,
	%\newblock {\em Fehlerfortpflanzung bei {I}nterpolation},
	%\newblock {Numer. Math.}, {\bf 3} (1961), pp. 62--71.
	
	\bibitem{Rivlin}
	T.~J. Rivlin,
	\newblock{\em Chebyshev polynomials},
	\newblock {Wiley}, New York, 2nd edition, 1990.
	
	\bibitem{Robidoux2006}
	N.~Robidoux,
	\newblock {\em Polynomial histopolation, superconvergent degrees of freedom, and pseudospectral discrete Hodge operators},
	\newblock {technical report} (2006). Available at \url{https://citeseerx.ist.psu.edu/document?repid=rep1&type=pdf&doi=7d9598cb27819bc85cedeb149368d4e213ce9f1c}

 	\bibitem{Runge}
	C. Runge,
	\newblock {\em Über empirische {F}unktionen und die {I}nterpolation zwischen äquidistanten {O}rdinaten},
	\newblock {Zeitschrift fur Mathematik und Physik}, {\bf 46}, (1901), pp. 224--243.
	
%	\bibitem{Schoenberg}
%	I. J. Schoenberg, 
%	\newblock {\em Splines and Histograms}, (1973)
%	\newblock in {Spline Functions and Approximation Theory}, A. Meir and A. Sharma, eds., ISNM International Series of Numerical Mathematics / Internationale Schriftenreihe zur Numerischen Mathematik / Série Internationale d’Analyse Numérique {\bf 21}, Birkhäuser, Basel, 2023, pp. 227--327.
%	
    \bibitem{Sommariva2009}
    A. Sommariva and M. Vianello,
    \newblock {Computing approximate {F}ekete points by {QR} factorizations of {V}andermonde matrices},
    \newblock {Comput. Math. Appl.}, {\bf 57}(8) (2009), pp. 1324--1336.
	
	\bibitem{Sundermann}
	B.~S\"{u}ndermann,
	\newblock{\em Lebesgue constants in {L}agrangian interpolation at the
		{F}ekete points},
	\newblock{Mitt. Math. Ges. Hamburg}, {\bf 11}(2) (1983), pp. 204--211.
	
	\bibitem{Szego}
	G. Szeg\H{o},
	\newblock{\em Orthogonal polynomials},
	\newblock {American Mathematical Society Colloquium Publications}, Providence, 1975.
	
	%\bibitem{trefappr}
	%L.~N. Trefethen, 
	%\newblock {\em Approximation Theory and Approximation Practice}, 
	%\newblock SIAM, Philadelphia, 2013.
	
	\bibitem{Trefethen}
	L.~N. Trefethen and J.~A.~C. Weideman,
	\newblock {\em Two results on polynomial interpolation in equally spaced points},
	\newblock {J. Approx. Theory}, {\bf 65}(3) (1991), pp. 247--260.
	
	%\bibitem{Vinogradov01}
	%O.~L. Vinogradov,
	%\newblock {\em Upper bounds of the Lebesgue constants for Fourier-Jacobi %series summation methods defined by a multiplier function},
	%\newblock {J. Math. Sci.}, {\bf 120}(5) (2004), pp. 1662--1671.
	
	%\bibitem{Whitneybook}
	%H. Whitney,
	%\newblock {\em Geometric integration theory},
	%\newblock {Princeton University Press, Princeton}, 1957.
	
	
	%	\bibitem{ZAR}
	%	E. Zampa, A. Alonso Rodr\'iguez and F. Rapetti,
	%	\newblock {\em Using the {FES} framework to derive new physical degrees of
		%		freedom for finite element spaces of differential forms},
	%	\newblock {Adv. Comput. Math.}, {\bf 49}(2):17 (2023).
	
	
\end{thebibliography}
\end{document}